\documentclass[aos,MSNbibl,seceqn,nameyear,dvips]{arximspdf}

\doi{10.1214/13-AOS1158} 
\volume{41}
\issue{6}
\pubyear{2013}
\firstpage{2739}
\lastpage{2767}

\makeatletter

\newtheorem{theo}{Theorem}
\newtheorem{Lem}{Lemma}
\newproclaim{rem}{Remark}

\newcommand{\cal}{\mathcal}

\def\log{\operatorname{log}}
\def\var{\operatorname{var}}

\def\a{\alpha}
\def\bDe{{\bar\De}}
\def\be{\beta}
\def\bg{{\bar g}}
\def\bow{\bowtie}

\def\cC{{\cal C}}
\def\cD{{\cal D}}
\def\cE{{\cal E}}
\def\cG{{\cal G}}
\def\cI{{\cal I}}
\def\cov{\operatorname{cov}}
\def\De{\Delta}
\def\ep{\varepsilon}
\def\err{\mathrm{err}}

\def\Ga{\Gamma}
\def\hal{{\widehat\a}}
\def\hbe{{\widehat\be}}
\def\hG{{\widehat G}}
\def\hGa{{\widehat\Ga}}
\def\hmu{{\widehat\mu}}
\def\hpsi{{\widehat\psi}}
\def\hth{{\widehat\th}}
\def\hsi{{\widehat\si}}
\def\ka{\kappa}
\def\kj{_{kj}}
\def\kji{_{kji}}
\def\kl{_{k\ell}}
\def\mi{\mid}
\def\om{\omega}
\def\ra{\to}
\def\rai{\ra\infty}
\def\scale{_{\mathrm{scale}}}
\def\si{\sigma}
\def\ssum{_{\mathrm{sum}}}

\def\th{\theta}
\def\tsi{{\widetilde\si}}

\makeatother

\begin{document}
\begin{frontmatter}

\title{Unexpected properties of bandwidth choice when smoothing
discrete data for constructing a~functional data classifier}
\runtitle{Bandwidth choice in classification}

\begin{aug}
\author[A]{\fnms{Raymond J.} \snm{Carroll}\thanksref{t3}\ead[label=e3]{carroll@stat.tamu.edu}},
\author[B]{\fnms{Aurore} \snm{Delaigle}\corref{}\thanksref{t2}\ead[label=e1]{A.Delaigle@ms.unimelb.edu.au}}
\and
\author[B]{\fnms{Peter} \snm{Hall}\thanksref{t2}\ead[label=e2]{halpstat@ms.unimelb.edu.au}}
\runauthor{R. J. Carroll, A. Delaigle and P. Hall}
\affiliation{Texas A\&M University, University of Melbourne and\break
University of Melbourne}
\address[A]{R. J. Carroll\\
Department of Statistics\\
Texas A\&M University\\
College Station, Texas 77843\\
USA\\
\printead{e3}}
\address[B]{A. Delaigle\\
P. Hall\\
Department of Mathematics and Statistics\\
University of Melbourne, Parkville\\
Victoria 3010\\
Australia\\
\printead{e1}\\
\hphantom{E-mail: }\printead*{e2}}
\end{aug}

\thankstext{t3}{Supported by a Grant from the National Cancer
Institute (R37-CA057030).
This publication is based in part on work supported by Award Number
KUS-CI-016-04, made by King Abdullah University of Science and
Technology (KAUST).}

\thankstext{t2}{Supported by grants and fellowships from the
Australian Research Council.}

\received{\smonth{3} \syear{2013}}
\revised{\smonth{7} \syear{2013}}

%
\begin{abstract}
The data functions that are studied in the course of functional data
analysis are assembled from discrete data, and the level of smoothing
that is used is generally that which is appropriate for accurate
approximation of the conceptually smooth functions that were not
actually observed. Existing literature shows that this approach is
effective, and even optimal, when using functional data methods for
prediction or hypothesis testing. However, in the present paper we
show that this approach is not effective in classification problems.
There a useful rule of thumb is that undersmoothing is often desirable,
but there are several surprising qualifications to that approach.
First, the effect of smoothing the training data can be more
significant than that of smoothing the new data set to be classified;
second, undersmoothing is not always the right approach, and in fact in
some cases using a relatively large bandwidth can be more effective;
and third, these perverse results are the consequence of very unusual
properties of error rates, expressed as functions of smoothing
parameters. For example, the orders of magnitude of optimal smoothing
parameter choices depend on the signs and sizes of terms in an
expansion of error rate, and those signs and sizes can vary
dramatically from one setting to another, even for the same classifier.
\end{abstract}

%

\end{frontmatter}

\section{Introduction}\label{sec1}
All supposedly ``functional'' data are actually observed discretely,
sometimes on a grid and on other occasions at randomly scattered
points. For example, in longitudinal data\vadjust{\goodbreak} analysis the observation
points are often widely spaced and irregularly placed, and substantial
smoothing is commonly used to convert discrete data like these to
functions. The impact of such smoothing has been addressed in the
context of prediction or hypothesis testing for functional data; see,
for example, \citet{HalVan07}, \citet{PanKraMad10},
\citet{WuMul11}, \citet{BenDeg}, \citet{CarJos11} and
\citet{CarDegJos}. The main conclusion of these papers has been
that conventional rules for smoothing discrete data typically apply,
and that smoothing parameters of standard size generally are
appropriate.

In contrast, the present paper was motivated by numerical work
indicating that, in the context of classifying functional data,
smoothing parameters of highly nonstandard sizes are appropriate, and
more generally that, even for a relatively simple classifier, there is
no simple precept (even an asymptotic prescription of size) that leads
to minimisation of error rate. If one had to give a rule, valid in some
but by no means all cases, it would be to undersmooth, but even there
unexpected caveats must be addressed.

For example, it turns out that the impact of smoothing the training
data can be more significant than that of smoothing the new data to be
classified. Indeed, the effect of smoothing the new data is
characteristic of a parametric problem, rather than a nonparametric
one. There, asymptotic arguments indicate that $(\mbox{sample
size})^{-1/2}$ is an appropriate bandwidth size for reducing the impact
of smoothing to parametric levels, whereas $(\mbox{sample
size})^{-1/3}$ is the nearest analogue for smoothing the training data.

However, both these recommendations are incorrect in many cases.
Depending on the signs and sizes of certain functionals of the data
distributions, it can be optimal to use smoothing parameters that are
an order of magnitude smaller, or an order of magnitude larger, than
these. Using some viewpoints the need for a low level of smoothing is
intuitively clear. Indeed, we expect that relatively minor features of
a curve, of the sort that might disappear if we were too enthusiastic
in the smoothing step, could have important information to convey in a
classification analysis. On the other hand, our results show that very
high levels of smoothing are sometimes advantageous.

We drew these conclusions after studying three different classifiers
for functional data: the standard centroid-based method, the
scale-adjusted form of that approach, and a version for functional data
of quadratic discrimination. Our conclusions are valid for all three
approaches, although they contradict conclusions which are well known
for standard nonparametric approximations to the Bayes classifier in
multivariate, rather than functional-data, settings. Specifically, for
univariate and functional data, and nonparametric Bayes classifiers,
conventional smoothing parameters, for example, those chosen using
standard plug-in rules for function estimation, typically are of the
correct order even though they do not quite minimise asymptotic
classification error; see, for example, \citet{HalKan05}. Moreover,
there does not exist a version of our results in univariate or
multivariate settings, since there is no analogue in such cases of the
``lattice effect,'' represented by the $m\kj$'s.

To these comments, we should add that in practice there is relatively
little difficulty in choosing smoothing parameters to minimise error
rate; cross-validation is usually effective. Our aim in this paper is
therefore not to develop methods for choosing the bandwidth optimally,
or nearly optimally, in classification problems, but to provide an
understanding of the many aspects of those problems that conspire
together to determine the optimal choice.

\section{Model and methodology}\label{sec2}

\subsection{Model}\label{sec21}

We consider $n_0$ (resp., $n_1$) unknown random functions $\{
g_{0j},\break 1\leq j\leq n_0\}$ (resp., $\{g_{1j},  1 \leq j \leq n_1\}
$) coming from two populations. We observe a training sample of the
data pairs $\cD\kj=\{(X\kji,Y\kji)$, $1\leq i\leq m\kj\}$, for
$1\leq j\leq n_k$ and $k = 0, 1$, corresponding to noisy versions of
the $g_{kj}$'s sampled at a discrete set of random points (i.e., $X\kji
$'s) and generated by the model
%
%
\begin{equation}\label{eq21}
Y\kji=g\kj(X\kji)+\ep\kji,
\end{equation}
where $k$ indexes the population, $\Pi_k$, from which the data in $\cD
\kj$ came, $j$ denotes the index of an individual drawn from $\Pi_k$,
and $i$ is the index of a data pair $(X\kji,Y\kji)$ for the $j$th
individual from the $k$th population.

The $g\kj$ are random functions defined on a compact interval $\cI$,
but observed only at $m\kj$ points $X_{kj1},\ldots,X_{kjm\kj}$.
These points may be fixed or random, and although we shall develop our
arguments in the random case, they can easily be extended to the fixed
case. We assume that each $g\kj$ has two bounded derivatives on $\cI
$; the respective sequences of $X$'s and $\ep$'s are each identically
distributed with distributions that do not depend on the $g$'s; the
$g$'s, $X$'s and $\ep$'s are all mutually independent; the $X$'s are
supported on $\cI$; and the $\ep$'s have zero mean and finite variance.

We also observe a new data set $\cD=\{(X_i,Y_i)$, $1\leq i\leq m\}$,
similar to the $\cD\kj$'s except that in this case we do not know
which population the data come from. Here,
%
%
\begin{equation}\label{eq22}
Y_i=g(X_i)+\ep_i,
\end{equation}
where the function $g$, the $X$'s and the $\ep$'s have the properties
given in the previous paragraph. Using the training data, we wish to
determine whether $\cD$ came from $\Pi_0$ or $\Pi_1$.

In the functional data literature [see, e.g., \citet{RamSil05}], when the data are noisy, it is common to preprocess them prior
to further analysis.
Typically, this is done by smoothing the data in some way, for example,
through a spline or kernel smoother, thereby obtaining, from the data
in $\cD\kj$ and $\cD$, estimators ${\widehat g}\kj$ and ${\widehat
g}$ of $g\kj$ and $g$, respectively. In the classification context,
once these estimators have been derived, they are plugged into
functional data classifiers, replacing there the unobserved functions
$g$ and $g\kj$ by their estimators ${\widehat g}$ and ${\widehat g}\kj
$. Our aim in this paper is to describe the application of estimators
${\widehat g}$ and ${\widehat g}\kj$ of $g$ and $g\kj$, and in
particular to describe the influence of tuning parameters used to
construct them, when the aim is classification rather than just
function estimation.

\subsection{\texorpdfstring{Estimating $g$, $g\kj$ and their mean and covariance
functions}{Estimating g, g kj and their mean and covariance
functions}}\label{sec22}

There are several ways to obtain nonparametric estimators of the
functions $g$ and $g\kj$, but the most popular ones are spline and
local linear methods. They have similar properties, but since local
linear estimators are much more tractable theoretically, we shall use
these in this work.
For $x\in\cI$, the local linear estimators of $g$ and $g\kj$ are
defined by
%
%
\begin{eqnarray}\label{eq23}
{\widehat g}(x)&=&{U_2(x) V_0(x)-U_1(x) V_1(x)\over U_2(x)
U_0(x)-U_1^2(x)},
\nonumber
\\
{\widehat g}\kj(x)&=&{U_{kj2}(x) V_{kj0}(x)-U_{kj1}(x)
V_{kj1}(x)\over U_{kj2}(x) U_{kj0}(x)-U_{kj1}^2(x)},
\end{eqnarray}
where
%
%
\begin{eqnarray}
\label{eq24}
U_\ell(x)&=&{1\over m} \sum_{i=1}^m
\biggl({x-X_i\over h} \biggr)^{
\ell} K_h
(x-X_i ),
\\
\label{eq26}
V_\ell(x)&=&{1\over m} \sum_{i=1}^m
Y_i \biggl({x-X_i\over h} \biggr)^{ \ell}
K_h ({x-X_i} ),
\\
\label{eq25}
U_{kj\ell}(x)&=&{1\over m\kj} \sum_{i=1}^{m\kj}
\biggl({x-X\kji\over h_1} \biggr)^{ \ell} K_{h_1} (x-X\kji ),
\\
V_{kj\ell}(x)&=&{1\over m\kj} \sum_{i=1}^{m\kj}
Y\kji \biggl({x-X\kji\over h_1} \biggr)^{ \ell} K_{h_1}
({x-X\kji} ),
\nonumber
\end{eqnarray}
$K$ is a kernel function, $h>0$ and $h_1>0$ are bandwidths, and
$K_h(x)=K(x/h)/h$. See, for example, \citet{FanGij96}. For
simplicity, throughout we use the same bandwidth $h_1$ for each
population and each individual, but we could have replaced $h_1$ by
bandwidths that depended on $k$ and $j$, as we do in our numerical work.

The classifiers we consider in this work require estimators of the
population means and covariances.
For $k=0,1$, let $\mu_k$ denote the mean function
%
%
\begin{equation}\label{eq27}
\mu_k=E_k(g)=E_k(g\kj),
\end{equation}
where $E_k$ represents expectation under the assumption that the data
come from $\Pi_k$. Also, let $G_k$ be the covariance function, defined by
$G_k(u,v)=\cov_k\{g(u),g(v)\}=E_k\{g(u) g(v)\} - \mu_k(u) \mu_k(v)$,
where $\cov_k$ denotes covariance when the data come from $\Pi_k$.
Estimators $\hmu_k$ and $\hG_k$ of $\mu_k$ and $G_k$ are defined in
the standard way by the empirical mean and covariance functions, but
replacing, in the definitions of these estimators, the unobserved $g\kj
$ by~${\widehat g}\kj$:
%
%
\begin{eqnarray}
\label{eq28}
\hmu_k&=&{1\over n_k} \sum_{j=1}^{n_k}
{\widehat g}\kj,
\\
\label{eq212}
\hG_k(u,v) &=& {1\over n_k} \sum
_{j=1}^{n_k} \bigl\{{\widehat g}\kj(u)-\hmu
_k(u) \bigr\} \bigl\{{\widehat g}\kj(v)-\hmu_k(v) \bigr\}.
\end{eqnarray}
See, for example, \citet{RamSil05}, Chapter~2.

Consider the spectral decomposition of the covariance function
%
%
\begin{equation}\label{eq215}
G_k(u,v)=\sum_{\ell= 1}^{\infty}
\th_{k\ell} \psi_{k\ell}(u) \psi_{k\ell}(v),
\end{equation}
where $(\th_{k\ell},\psi_{k\ell})$ is an (eigenvalue,
eigenfunction) pair for the linear operator $G_k$ defined by $G_k(\psi
)(u)=\int G_k(u,v) \psi(v) \,dv$, and where, following convention, we
have used the notation $G_k$ for both the operator and the covariance.
The terms in (\ref{eq215}) are ordered such that $\th_{k1}\geq\th
_{k2}\geq\cdots\geq0$.
If $g$ is drawn from $\Pi_k$ then we can write
%
%
\begin{equation}\label{eq214}
g(x)=\mu_k(x)+\sum_{\ell= 1}^{\infty}
Z_{k\ell} \th_{k\ell
}^{1/2}\psi_{k\ell}(x),
\end{equation}
where $\mu_k=E_k(g)$ denotes the mean of the random process of which
$g$ is a realisation, $Z_{k\ell}=\th_{k\ell}^{-1/2}\int(g-\mu_k)
\psi_{k\ell}$, and the $Z_{k\ell}$'s (for $\ell=1,2,\ldots$)
comprise a sequence of uncorrelated random variables with zero mean and
unit variance.
The quantities\vspace*{1pt} $\th_{k\ell}$ and $\psi_{k\ell}$ can be estimated
consistently by the eigenvalues and eigenfunctions $\hth_{k\ell}$ and
$\hpsi_{k\ell}$ of the linear\vspace*{1pt} operator $\hG_k$, defined by $\hG
_k(\psi)(u)=\int\hG_k(u,v) \psi(v) \,dv$, with the covariance
estimator $\hG_k$ defined as at (\ref{eq212}):
%
%
\begin{equation}\label{eq217}
\hG_k(u,v) =\sum_{\ell=1}^{\infty}
\hth_{k\ell} \hpsi_{k\ell}(u) \hpsi _{k\ell}(v),
\end{equation}
where $\hth_{k1}\geq\hth_{k2}\geq\cdots\geq0$, and, since $\hth
_{k\ell}=0$ for all $\ell>n_k$, all but the first $n_k$ terms in the
series at (\ref{eq217}) vanish. See \citeauthor{HalHos06}
(\citeyear{HalHos06,HalHos09})
for properties of these estimators in the case where $g$ and $g\kj$
are observed; see also Li and Hsing (\citeyear{LiHsi10N1,LiHsi10N2}) for other cases.

\subsection{Constructing classifiers}
Classifiers for functional data have received a great deal of attention
in the literature. See, for example, \citet{VilPer04},
\citet{BiaBunWeg05}, \citet{FroTul06}, \citet{LenMul06},
\citet{LopRom06}, \citet{RosVil06}, \citet{CueFebFra07},
\citet{WanRayMal07}, \citet{BerBiaRou08}, \citet{Epi08},
\citet{Araetal09}, \citet{DelHal12} and \citet{DelHalBat12}.

In those papers the authors suggest methods for constructing
classifiers, but so far the theoretical impact of smoothing; that is,
the impact of using ${\widehat g}$ and ${\widehat g}\kj$ instead of
$g$ and $g\kj$ when constructing classifiers; has been largely ignored
in the literature. In this paper, we study this impact of smoothing for
three relatively simple functional classifiers: the centroid
classifier, or Rocchio classifier [see, e.g., \citet{ManRagSch08}],
commonly used for classifying high-dimensional data; a scaled version
of this classifier, which we define below in a general way; and a
version for functional data of Fisher's quadratic discriminant,
studied,
for example, by \citet{LenMul06} and \citet{DelHal12}.
These classifiers are usually defined in terms of the functions $g$ and
$g\kj$, and here we shall define them in terms of ${\widehat g}$ and
${\widehat g}\kj$. The standard versions of these classifiers are
obtained by replacing ${\widehat g}$ and ${\widehat g}\kj$ by $g$ and
$g\kj$. The functions ${\widehat g}\kj$ appear only implicitly
through the estimated means and covariance functions constructed in
Section~\ref{sec22}.

In the present setting, the centroid-based classifier assigns the curve
$g$, observed through $\cD$, to $\Pi_0$ if the statistic
%
%
\begin{equation}\label{eq29}
S({\widehat g})=\int_\cI \bigl\{{\widehat g}(t)-
\hmu_0(t)\bigr\}^2 \,dt -\int_\cI
\bigl\{{\widehat g}(t)-\hmu_1(t)\bigr\}^2 \,dt
\end{equation}
is negative, and to $\Pi_1$ if $S({\widehat g})>0$.

A scaled version of the centroid classifier, which accommodates
differences in scales between the two populations, can be defined by
replacing $S$ in (\ref{eq29}) by
%
%
\begin{eqnarray}\label{eq210}
S\scale({\widehat g}) &=&{1\over s_0^2} \int_\cI
\bigl\{{\widehat g}(t)-\hmu_0(t)\bigr\}^2 \,dt \nonumber\\[-8pt]\\[-8pt]
&&{}-
{1\over s_1^2} \int_\cI \bigl\{{\widehat g}(t)-
\hmu_1(t)\bigr\}^2 \,dt +\log \biggl({s_0^2\over s_1^2}
\biggr),\nonumber
\end{eqnarray}
where $s_k^2$ is an estimator of the scale of population $\Pi_k$. For
example, we might take $s_k^2$ to equal
${n_k}^{-1}\sum_{j=1}^{n_k}
\int_\cI ({\widehat g}\kj-\hmu_k)^2$,
the version we used in our numerical work, or
$\int_\cI\int_\cI\hG_k(u,v) \psi(u) \psi(v) \,du \,dv, %
$
where
$\psi$ is open to choice; or $s_0^2$ and $s_1^2$ could be selected
empirically by minimising a cross-validation estimator of
classification error. The definition at (\ref{eq210}) should be
compared with those at (\ref{eq213}) and (\ref{eq216}), below. The
form of (\ref{eq210}), and also of (\ref{eq213}) and (\ref
{eq216}), is motivated by likelihood-ratio statistics for Gaussian data.

A version for functional data of Fisher's quadratic discriminant is
based on
%
%
\begin{eqnarray}\label{eq213}
T({\widehat g})&=&\sum_{\ell= 1}^{p} \biggl[
{1\over\hth_{0\ell}} \biggl\{\int_\cI ({\widehat g}-
\hmu_0) \hpsi_{0\ell} \biggr\}^{ 2} \nonumber\\[-8pt]\\[-8pt]
&&\hspace*{19.2pt}{}-
{1\over\hth_{1\ell}} \biggl\{\int_\cI ({\widehat g}-
\hmu_1) \hpsi _{1\ell} \biggr\}^{ 2} +\log \biggl(
{\hth_{0\ell}\over\hth_{1\ell}} \biggr) \biggr],\nonumber
\end{eqnarray}
where ${\widehat g}$ and $\hmu_k$ are as at (\ref{eq23}) and (\ref
{eq28}), $(\hth_{0\ell},\hth_{1\ell})$ are at (\ref{eq217}) and
$p$ is a positive truncation parameter. (Here we assume, as is often
the case in practice, that the prior probabilities of each population
are unknown and estimated by $1/2$. A more general version of the
classifier can be used if these probabilities are estimated by other
values, but this does not alter our main conclusions.) We assign the
new data set $\cD$ to $\Pi_0$ if $T({\widehat g})\leq0$, and to $\Pi
_1$ otherwise.
Of course, the statistic $T({\widehat g})$, at~(\ref{eq213}), is just
an empirical version of the quantity
%
%
\begin{eqnarray}\label{eq216}
T_0(g) &=&\sum_{\ell= 1}^{p} \biggl[
{1\over\th_{0\ell}} \biggl\{\int_\cI (g-
\mu_0) \psi_{0\ell} \biggr\}^{ 2}\nonumber\\[-8pt]\\[-8pt]
&&\hspace*{18.2pt}{} -
{1\over\th_{1\ell}} \biggl\{\int_\cI (g-
\mu_1) \psi_{1\ell
} \biggr\}^{ 2} +\log \biggl(
{\th_{0\ell}\over\th_{1\ell}} \biggr) \biggr].\nonumber
\end{eqnarray}

If the functions $g$ are Gaussian, and the first $p$ eigenvalues, in
versions of (\ref{eq215}) and (\ref{eq217}) for either population,
are distinct and nonzero, and the remaining eigenvalues vanish, then
the classifier based on $T_0(g)$, at (\ref{eq216}), is optimal in the
sense of having least classification error among all classifiers, since
it is, after all, just a likelihood ratio statistic. When the
eigenvalues and eigenfunctions are estimated from data, as at (\ref
{eq213}), the classifier is asymptotically optimal. Bearing in mind the
effectiveness of Fisher's discriminant analysis in the case of
vector-valued data, even when the data are not normal, the classifier
based on $T({\widehat g})$ is an attractive choice even in non-Gaussian cases.

\section{Theoretical properties}\label{sec3}

\subsection{Standard centroid-based classifier}\label{sec31}
In this section, we derive properties of the centroid classifier based
on the estimators ${\widehat g}$ and ${\widehat g}\kj$, and in
particular we examine the impact of smoothing.
First, we introduce notation.
Let $n=n_0+n_1$ (hence $n$ is a positive integer sequence diverging to
infinity), let $m=m(n)$ be of the same size as $m\kj$ [see (\ref
{eq39}) below], and write $\si_{\ep k}^2$ for the variance of
the\vadjust{\goodbreak}
experimental errors $\ep\kji$ and $\ep_i$, in (\ref{eq21}) and
(\ref{eq22}), when the data come from $\Pi_k$.
Let
%
%
\begin{equation}\label{eq31}
\bg_k={1\over n_k} \sum_{j=1}^{n_k}
g\kj,\qquad \nu_k=n_k^2 \Biggl(\sum
_{j=1}^{n_k} m\kj^{-1} \Biggr)^{ -1}
\end{equation}
and define
%
%
\begin{eqnarray}
\label{eq33}
b_{k0}&=&\int_\cI (\mu_1-
\mu_0) \bigl\{2 \mu_k-(\mu_0+
\mu_1)\bigr\},\nonumber\\[-8pt]\\[-8pt]
\be_{k0}&=&\int_\cI (
\bg_1-\bg_0) \bigl\{2 \mu_k-(
\bg_0+\bg_1)\bigr\},\nonumber
\\
\label{eq35}
\si_k^2&=&4 \ka_2\int_\cI
\int_\cI (\bg_1-\bg_0)
(x_1) (\bg_1-\bg_0) (x_2)
G_k(x_1,x_2) \,dx_1
\,dx_2,
\\
\label{eq36}
\tau_k^2&=&4 \ka_2\int_\cI
\int_\cI (\mu_1-\mu_0)
(x_1) (\mu_1-\mu_0) (x_2)
G_k(x_1,x_2) \,dx_1
\,dx_2,
\end{eqnarray}
where $\ka_2=\int u^2 K(u) \,du$. Finally, put $\ka=\int K^2$, and
let $\cI$ be the compact interval that equals the support of the
density $f_X$ of the $X_i$'s and $X\kji$'s, and of the functions $g$
and $g\kj$.

We make the following assumptions:
%
%
\begin{eqnarray}\label{eq37}
\begin{tabular}{p{320pt}} (a) The distribution of the
experimental errors $\ep\kji$ and $\ep _i$, in (\ref{eq21})
and~(\ref{eq22}), has zero mean and all moments finite, may depend on
$k$, and has variance $\si_{\ep k}^2$; (b) the density $f_X$ of the
variables $X\kji$ and $X_i$ does not depend on $i$, $j$ or $k$; (c)
$f_X$ has two bounded derivatives, $f_X(x)\geq C>0$ for all $x\in\cI$,
and $f_X''$ is H\"older continuous on the support $\cI$ of $f_X$.
\end{tabular}\hspace*{-34pt}
\\
%
\label{eq38}
\begin{tabular}{p{320pt}} (a) The functions $g$ and $g\kj$
associated with the populations $\Pi _k$, for $k=0,1$, are
realisations of Gaussian processes, have uniformly bounded covariance functions
$G_k$ and mean functions $\mu_k$, both depending only on
$k$, and satisfy $\tau_k^2>0$ for $k=0,1$; and (b) with
probability 1 the functions are uniformly bounded and have H\"older-continuous
second derivatives, with the property that, for a constant $C>0$, all moments of
$\sup_{x_1,x_2} |g''(x_1)-g''(x_2)|/|x_1-x_2|^C$
are finite\vspace*{1pt} when $g$ is sampled from either $\Pi_0$ or $
\Pi_1$.
\end{tabular}\hspace*{-34pt}
\\
%
\label{eq39}
\begin{tabular}{p{320pt}} (a) For a constant $C>0$, the
results $h^{(1)}=O(n^{-C})$ and $n^{1-C}
h^{(1)}\rai$ hold for $h^{(1)}=h$ and $h^{(1)}=h_1$;
(b) the kernel $K$ is a symmetric, nonnegative, compactly supported and H\"older
continuous probability density; and (c) for each $k$ the values of
$m^{-1}\min_jm\kj$, $m^{-1}\max
_jm\kj$ and $n_0/n_1$ are bounded
away from zero and infinity as $n\rai$, and, for constants $C_1$ and
$C_2$ satisfying $0<C_1<$ $C_2<\infty$, $m$
and $n_0$ lie between $n^{C_1}$ and~$n^{C_2}$.
\end{tabular}\hspace*{-34pt}
\end{eqnarray}

The assumption in (\ref{eq37})(c) that $f_X$ is bounded away from zero on its
support is only a technical requirement, and is unnecessary in
practice. To make this clear, in our numerical work we shall take $f_X$
to be a normal density, and show that the conclusions of Theorem
\ref{Theorem1} are nevertheless reflected clearly.

Let $\err_k=P_k\{(-1)^kS({\widehat g}) >0\}$ denote the probability
that the standard centroid-based classifier, based on the statistic
$S({\widehat g})$ at (\ref{eq29}), commits an error when the data set
$\cD$ actually comes from $\Pi_k$.
Theorem \ref{Theorem1} below describes the asymptotic behaviour of
$\err_k$, and highlights the effect of the smoothing parameters $h$
and $h_1$, used to construct the estimators ${\widehat g}$ and
${\widehat g}\kj$ of $g$ and $g\kj$, on the classifier. A proof is
given in Appendix \ref{app11}.

\begin{theo}\label{Theorem1}
Assume that (\ref{eq37})--(\ref{eq39}) hold, and let $\Psi
_0=1-\Phi$ and $\Psi_1=\Phi$, where $\Phi$ denotes the c.d.f. of a
standard normal random variable. Then
%
%
\begin{eqnarray}\label{eq310}
\err_k &=& \err_{k0} + h^2 c_k
+h_1^2 c_{k1} +\frac{d_{k0}}{\nu_0h_1} +
\frac{d_{k1}}{\nu_1h_1}
\nonumber\\[-8pt]\\[-8pt]
&&{}+O \bigl\{m^{-1}+(mh)^{-2} \bigr\} +o
\biggl(h^2+h_1^2+\frac{1}{\nu_0h_1} \biggr),\nonumber
\end{eqnarray}
where
$\err_{k0}=E_k[\Psi_k\{-\be_{k0}/\si_k\}]$,
$c_k=\ka_2 \alpha_k\int_\cI (\mu_1-\mu_0) \mu_k''$,
$c_{k1}=-\ka_2 \alpha_k \allowbreak\int_\cI (\mu_1-\mu_0) \mu_{1-k}''$,
$d_{kj}=(-1)^{j} \alpha_k \si_{\ep j}^2 \ka\int_\cI f_X^{-1}$,
with $\alpha_k=(-1)^k \tau_k^{-1}\phi(b_{k0}/\tau_k)$,
and where $\phi$ denotes the standard normal density function.
\end{theo}

The leading term $\err_{k0}$ on the right-hand side of (\ref{eq310})
does not depend in any way on the bandwidths $h$ and $h_1$. It does
involve the training sample sizes $n_0$ and~$n_1$, and in particular
does not equal the asymptotic limit of $\err_k$ as $n$ increases,
since that limit is given by $\Psi_k(-b_{k0}/\tau_k)$, but the
effects of the bandwidths are all confined to subsequent terms on the
right-hand side of (\ref{eq310}).
The terms in $h^2$ and $h_1^2$ represent contributions to
classification error arising from biases of the estimators ${\widehat
g}$ and ${\widehat g}\kj$,
and the terms in $(\nu_0h_1)^{-1}$ and $(\nu_1h_1)^{-1}$ are
contributions from the variances of the estimators ${\widehat g}\kj$.

While a priori it might be thought that, since the total number of
observations in the training sample, $\sum_j m\kj$, for $k=0$ and 1, is
an order of magnitude larger than the number of observations, $m$, in
the new data set $\cD$, then $h_1$ should be chosen smaller than $h$,
Theorem \ref{Theorem1} shows that the influence of bandwidths on error
rate is much more complex than this.

For one thing, there are no terms in $(mh)^{-1}$ on the right-hand side
of (\ref{eq310}). (Section~\ref{Remark4} will explain the reason for
this.) As a result, the terms on the right-hand side of (\ref{eq310})
that depend on $h$ can be rendered equal to $O(m^{-1})$ simply by taking
$h$ equal to a constant multiple of $m^{-1/2}$. As noted in Remark \ref
{Remark7}, below, this level of contribution to the error rate is
generally impossible to remove, even in simple parametric
problems.\vadjust{\goodbreak}
Therefore the contribution of $h$ to error rate cannot be rendered
smaller than $m^{-1}$. However, in some instances choosing $h$ to be an
order of magnitude larger or smaller than $m^{-1/2}$ can be beneficial;
see Section~\ref{Remark2} below.

The terms in $h_1$ on the right-hand side of (\ref{eq310}) are a
different matter because each of $c_{k1}$ and $d_{k0} \nu_0^{-1}
+d_{k1} \nu_1^{-1}$ can be either positive or negative. Depending on
the signs and sizes of $c_{k1}$ and $d_{k0} \nu_0^{-1}+d_{k1} \nu
_1^{-1}$, it can be optimal to take $h_1$ to be of order $\nu
_k^{-1/3}$, which achieves a trade-off between terms in $h_1^2$ and
$(\nu_kh_1)^{-1}$, or to take $h_1$ to decrease to zero more quickly or
to converge to a positive constant, as $n$ increases; see Section~\ref{Remark2} below.

Therefore, the impact that smoothing has on classification performance
is much more subtle than it might have appeared. We discuss these
issues in more detail in the next sections.

\subsubsection{Sizes of $h$ and $h_1$ that optimise overall error
rate}\label{Remark2}
Using Theorem \ref{Theorem1} we can deduce the orders of magnitudes of
$h$ and $h_1$ that minimise the error rate of the classifier, that is,
that minimise the probability of misclassification,
%
%
\begin{equation}\label{eq312}
\err=\pi_0 \err_0+\pi_1 \err_1,
\end{equation}
where $\err_0$ and $\err_1$ are as in (\ref{eq310}), $\pi_k$
denotes the prior probability attached to population $\Pi_k$, and $\pi
_0+\pi_1=1$.
Using (\ref{eq310}) and (\ref{eq312}), we can write
%
%
\begin{eqnarray}\label{eq312a}
\err&=&\err^0 +c^0 h^2+c_1^0
h_1^2 +d^0 (\nu_0h_1)^{-1}
+O \bigl\{m^{-1}+(mh)^{-2} \bigr\}
\nonumber\\[-8pt]\\[-8pt]
&&{}+o \bigl\{h^2+h_1^2+(\nu_0h_1)^{-1}
\bigr\},\nonumber
\end{eqnarray}
where $\err^0= \pi_0 \err_{00}+\pi_1 \err_{10}$ (recall that
$\err_{k0}$ does not depend on the bandwidths),
\begin{eqnarray*}
c^0&=& \ka_2\int(\mu_1-\mu_0)
\biggl\{\pi_0 {\mu_0''\over\tau_0} \phi \biggl(
{b_{00}\over\tau
_0} \biggr) -\pi_1 {\mu_1''\over\tau_1} \phi
\biggl({b_{10}\over\tau
_1} \biggr) \biggr\},
\\
c_1^0&=& \ka_2\int(\mu_0-
\mu_1) \biggl\{\pi_0 {\mu_1''\over\tau_0} \phi
\biggl({b_{00}\over\tau
_0} \biggr) -\pi_1 {\mu_0''\over\tau_1}
\phi \biggl({b_{10}\over\tau
_1} \biggr) \biggr\},
\\
d^0&=& \ka \biggl(\int_\cI
f_X^{-1} \biggr) \biggl\{{\pi_0\over\tau_0} \phi
\biggl({b_{00}\over\tau_0} \biggr) -{\pi_1\over\tau_1} \phi \biggl(
{b_{10}\over\tau_1} \biggr) \biggr\} \biggl(\si_{\ep0}^2-
\si_{\ep1}^2 {\nu_0\over\nu_1} \biggr).
\end{eqnarray*}
Since the function $\phi$ is symmetric, and $b_{10}=-b_{00}$, then
$b_{10}$ can be replaced by $b_{00}$ in the formula for $d_0$ without
altering its veracity.

To appreciate the very wide range of optimal bandwidth choices that can
arise in the problem of minimising error rate, let us consider
minimising $\err$, at (\ref{eq312a}).
To help remove ambiguities, let us assume that as $n$ increases the
value of $\si_{\ep0}^2-\si_{\ep1}^2 \nu_0 \nu_1^{-1}$ is of the
same sign for all sufficiently large $n$, and its absolute value is
bounded away from zero; assumption (\ref{eq39})(c) ensures that it is
uniformly bounded. In this instance, and focusing just on the terms in
$h_1$, we see that four distinct cases can arise in practice:
\begin{longlist}
\item $c_1^0$ and $d^0$ are both positive. In this case, to
minimise the contribution from~$h_1$, we should minimise $c_1^0
h_1^2+d^0 (\nu_0h_1)^{-1}$, which is achieved by taking $h_1$ to be of
size $\nu_0^{-1/3}$.

\item $c_1^0$ and $d^0$ are both negative. In this case, the
contribution made by $h_1$ behaves like $-\{|c_1^0| h_1^2+|d^0|(\nu
_0h_1)^{-1}\}$ as sample size increases. The term within braces here is
maximised by taking $h_1=0$, and analogously, in minimising $\err$, it
is optimal to take $h_1$ to be of strictly smaller order than~$\nu_0^{-1/3}$.

\item $c_1^0>0$ and $d^0\leq0$. In this case, to minimise the
error rate, we need to maximise the size of the negative term and
minimise that of the positive term, which is achieved by taking $h_1$
to be of strictly smaller order than $\nu_0^{-1/3}$ (the precise order
depends on the magnitude of second order terms, but deriving the latter
precisely would require a lot of additional computation).

\item $c_1^0<0$ and $d^0\geq0$. Here, using arguments similar
to those in case (iii), taking $h_1$ to be of strictly larger order
than $\nu_0^{-1/3}$ is optimal.
\end{longlist}
The case $d^0=0$ occurs, for example, if the covariance $G_k$ of the
Gaussian process~$g$, the experimental error variance $\si_{\ep k}^2$,
and the values of $m\kj$ and $n_k$ do not depend on $k$. Equal values
of $m\kj$ commonly arise when the data are observed on a grid; see
Remark \ref{Remark9}.

A similar analysis can be carried out in the case of optimisation over
$h$ rather than $h_1$, although there the optimum is accessed from a
comparison of terms in $h$ and $(mh)^{-2}$, rather than $h_1^2$
and $(\nu_0h_1)^{-1}$. [A tedious analysis of the term of size $(mh)^{-2}
$, represented by the remainder $O\{(mh)^{-2}\}$ in (\ref{eq310}),
shows that it can be either positive or negative.]
Depending on the relative signs of the terms in $h^2$ and $(mh)^{-2}$,
it can be optimal to take $h\asymp m^{-1/2}$, or $h$ of strictly larger,
or strictly smaller order than $m^{-1/2}$.

Similar results are obtained if we investigate properties of $\err_k$,
in (\ref{eq310}), instead of the overall error rate, $\err$,
at (\ref{eq312}).

These results explain the very diverse patterns of behaviour that are
seen in numerical work, and that motivated our research; see
Section~\ref{sec1}. In summary, in apparently similar problems and
using the same type of classifier, it can be optimal to use a very
small bandwidth, or a very large bandwidth, or a bandwidth of only
moderate size, depending on the signs of certain constants. Therein
lies the contradictory nature of the smoothing parameter choice problem
for classification of functional data.

\subsubsection{Absence of terms in $(mh)^{-1}$}\label{Remark4}

The centroid-based classifier statistic $S({\widehat g})$, at (\ref
{eq29}), can be written equivalently as
%
%
\begin{equation}\label{eq313}
S({\widehat g})=\int_\cI (\hmu_1-
\hmu_0) (2 {\widehat g}-\hmu_0-\hmu_1) \,dt.
\end{equation}
Importantly, there is no quadratic term in ${\widehat g}^2$ in (\ref
{eq313}), and as a result the impact of the bandwidth $h$, although not
$h_1$, on properties of the classifier is greatly reduced. This
reduction is brought about by the smoothing effect of the integral in
(\ref{eq313}), which results in the elimination of terms in $(mh)^{-1}$.

This property, to which we shall refer to as the ``integration
effect,'' is known in other settings, for example, when integrating a
kernel density estimator, computed from a sample of size $m$, to
produce a distribution estimator. Integration results in the variance
reducing from order $(mh)^{-1}$, for the density estimator, to order
$m^{-1}$, for the distribution function estimator---just as it does in
the setting above.

\begin{rem}[(Order $m^{-1}$ term in expansion of classification error)]\label{Remark7}
We assumed in (\ref{eq39})(c) that the values of $m\kj$,
representing the number of pairs $(X\kji,Y\kji)$ for a given
population index $k$ and given individual $j$, are all of roughly the
same size. In this setting it is easy to see that, even in an
elementary parametric setting, we must expect the operation of
observing the functions $g\kj$ at scattered points to affect error
rate through a term of order $m^{-1}$, and no smaller. For example,
consider the case where $g\kj=\psi( \cdot\mi\om\kj)$, with $\psi
( \cdot\mi\om)$ being a known function completely determined by the
parameter $\om$, and $\om\kj=\int_\cI g\kj w$ where the weight
function $w$ is known. Using the data $\cD\kj$ on $g\kj$ we can
estimate $\om\kj$ root-$m$ consistently, but no faster, and as a
result we incur a classification error of size $m^{-1}$, and no smaller,
from not knowing the values $\om\kj$. It is for this reason that,
when developing expansions of classification error, we do not explore
the remainder of size $m^{-1}$; it is stated simply as $O(m^{-1})$ on the
right-hand side of (\ref{eq310}).
\end{rem}

\subsubsection{Other remarks}\label{secother}
We conclude our discussion of Theorem \ref{Theorem1} with a number of remarks.

\begin{rem}[(Definition of $\hmu_k$)]\label{Remark5}
The size of the fourth and fifth terms on the right-hand side of (\ref
{eq310}) is determined by the sizes of $\nu_0^{-1}$ and $\nu_1^{-1}$,
and those quantities can be made slightly smaller by using a slightly
different definition of $\hmu_k$, at~(\ref{eq28}). In particular in
(\ref{eq28}), on account of the definition of ${\widehat g}\kj$ at
(\ref{eq23}), $\hmu_k$ is defined as an average of ratios of sums,
whereas slightly better statistical performance is obtained by taking
$\hmu_k$ to be simply a ratio of sums:
\[
\hmu_k={\sum_j (U_{kj2} V_{kj0}-U_{kj1} V_{kj1})\over
\sum_j (U_{kj2} U_{kj0}-U_{kj1}^2)},
\]
compare (\ref{eq23}). However, this approach departs from standard
practice in working with functional data, and therefore, since
convergence rates do not alter (only the constant multiples of rates
are reduced), we have followed standard practice in the definition
of $\hmu_k$.
\end{rem}

\begin{rem}[(Gaussian assumption)]\label{Remark8}
Of course, if $m$ is sufficiently large then ${\widehat g}$ is itself
approximately Gaussian, and so the assumption that $g$ is a Gaussian
process is reflected particularly well in properties of its estimator.
More generally, our assumption that $g$ is a Gaussian process is made
for simplicity, and can be relaxed. For example, generalisations to
chi-squared and other processes, where shape can be described in terms
of a small number of fixed functions (mean and covariance in the
Gaussian case), are straightforward.

More generally we would require a model which described the properties
of random functions relatively simply. The Gaussian model fills this
need ideally; shape is described by mean and variance functions, on
which we have imposed only smoothness, rather than parametric,
conditions. Moreover, in the Gaussian case all moments of $g(x)$ are
finite, for each $x$ (we use this property repeatedly during our
theoretical arguments), and the principal component scores are
independent (this is used frequently during our proof of Theorem \ref{Theorem2}).
\end{rem}

\begin{rem}[(Case of regularly spaced design)]\label{Remark9}
Theorem \ref{Theorem1} continues to hold if the $m\kj$ design
variables $X\kji$ are regularly spaced on $\cI$ for each $k$ and $j$.
The only change necessary is to replace $\int_\cI f_X^{-1}$, on the
right-hand side of (\ref{eq310}), by the square of the length of the
interval $\cI$.
\end{rem}

\subsection{Scale-adjusted centroid-based classifier}\label{sec32}
Recall that scale-adjusted centroid-based classifier is defined in
terms of $S\scale({\widehat g})$, at (\ref{eq210}).
A decomposition similar to that of Theorem \ref{Theorem1} can be
derived for this classifier, as we shall prove in Theorem \ref
{Theorem2} below.
For this classifier, it seems necessary to strengthen (\ref{eq39}) by
imposing conditions on the behaviour of the eigenvalues $\th\kl$ as
$\ell$ increases. However, since our aim in this section is only to
corroborate the conclusions in Section~\ref{sec31}, drawn there in
the case of the standard centroid-based classifier, then we shall
simplify our account by assuming that $g$ is finite dimensional, and in
particular taking the covariance expansion at (\ref{eq215}) to have
just $q$ terms:
%
%
\begin{equation}\label{eq314}
\begin{tabular}{p{320pt}} For $k=0$ and 1: (a) the first $q$
eigenvalues in the sequence $\th _{k1}\geq\th_{k2}\cdots\,$, arising in
the covariance expansion (\ref{eq215}) of $g$ when the data come
from~$\Pi_k$, are distinct; (b) $\th\kl=0$ for $\ell>q$; (c) for $1\leq
\ell\leq q$ the eigenfunctions $\psi\kl$ in (\ref{eq215}) have two
H\"older continuous derivatives on $\cI$; (d) $E_k(g)$ is a linear form
in $ \psi_{k1},\ldots,\psi_{kq}$; and (e) in the definition of $S
\scale({\widehat g})$, $s_0^2\neq s_1^2$.
\end{tabular}\hspace*{-34pt}
\end{equation}
Without (\ref{eq314})(a), separate conditions, valid uniformly in
$j=1,2,\ldots\,$, have to be imposed on remainders in Taylor expansions
of ``smoothed'' versions of the eigenvalues $\th\kj$, depending on $h$.

The next theorem indicates that the results of Theorem \ref{Theorem1}
also apply for the scale-adjusted centroid-based classifier. Its proof
is given in the supplementary material [\citet{CarDelHal}].

\begin{theo}\label{Theorem2} Assume that (\ref{eq37}), (\ref
{eq38}) and (\ref{eq314}) hold. Then the error rate of the
scale-adjusted centroid-based classifier, when the data in $\cD$ are
drawn from~$\Pi_k$, admits the expansion at (\ref{eq310}), but with
different constants, where the various terms have the properties stated
immediately below that formula.
\end{theo}

The diversity of possible signs of $c_k$, $c_{k1}$ and $d_{k0} \nu
_0^{-1}+d_{k1} \nu_1^{-1}$ in (\ref{eq310}), discussed in Section~\ref{Remark2}, is also present in this case.
Therefore the conclusions drawn in that section apply to the
scale-adjusted centroid-based classifier.
However, we have not derived explicitly the counterparts of the
constants $c_k$, $c_{k1}$, $d_{k0}$ and $d_{k1}$ that appear in
equation (\ref{eq310}).

The integration effect discussed in Section~\ref{Remark4} is also
present here, although we had originally expected that the
scale-adjusted centroid classifier would produce a term of size
$(mh)^{-1}$ in an expansion of error rate. Indeed, the situation
initially seems quite different in the case of the scale-adjusted
version $S\scale({\widehat g})$ of~$S({\widehat g})$, at (\ref
{eq210}), when $s_0^2\neq s_1^2$. There the quadratic term in
${\widehat g}$ persists. The reason it still does not produce a term in
$(mh)^{-1}$ is quite subtle. Define $\bow_k$ to be $>$ or $\leq$
according as $k=0$ or $k=1$, respectively. The probability $P_k\{
S\scale({\widehat g})\bow_k0\}$ can be written as
\[
P_k\Biggl\{\sum_{j=1}^\infty
w_j (Z_j+V_j)^2\bow_k
W\Biggr\} +\mbox{negligible terms},
\]
where the $Z_j$'s are independent $\mathrm{N}(0,1)$ variables, conditional on
the $V_j$'s and $W$; the positive weights $w_j$ are nonrandom; and
critically, $W$ does not involve the experimental errors $\ep_i$ in
(\ref{eq22}), from which any term in $(mh)^{-1}$ would arise. The
terms $V_j$ depend on the experimental errors only through integrals of
the error process, and the integration effect at this point largely
removes the impact of the error bandwidth $h$, with the result that
there is no term of size $(mh)^{-1}$. However, terms in $(\nu_0h_1)^{-1}
$ remain; the integration effect only influences smoothing of the new
data, not of the training data.

\subsection{Quadratic discriminant}\label{sec33}

Finally, we show that similar smoothing effects are present in the case
of the quadratic discriminant classifier defined through the statistic
$T({\widehat g})$ at (\ref{eq213}). Recall that, when the data in
$\cD$ come from $\Pi_k$, the random function $g$ has covariance
function $G_k$. To derive the counterpart of Theorem \ref{Theorem1} for
this classifier, let $r,r_1,r_2$ take the values 0 and 1, let
$1\leq\ell,\ell_{1},\ell_{2}\leq p$, and define the covariances
\[
\cov_k[r_1,r_2;\ell_{1},
\ell_{2}] =\int_\cI \int_\cI
G_k(x_1,x_2) \psi_{r_1\ell_{1}}(x_1)
\psi_{r_2\ell_{2}}(x_2) \,dx_1 \,dx_2,
\]
the variances $\var_k[r,\ell]=\cov_k[r,r;\ell,\ell]$, and the correlations
\[
\rho_k[r_1,r_2;\ell_{1},
\ell_{2}] ={\cov[r_1,r_2;\ell_{1},\ell_{2}]
\over(\var_k[r_1,\ell_{1}] \var_k[r_2,\ell_{2}])^{1/2}}.
\]

Let $p\geq1$, a fixed number, be the number of principal components
used to construct the quadratic discriminant statistic $T({\widehat
g})$, defined at (\ref{eq213}). Theorem \ref{Theorem3} below
addresses the error rate of the quadratic discriminant based on
$T({\widehat g})$, and there we shall assume that:
%
%
\begin{equation}\label{eq315}
\begin{tabular}{p{320pt}} (a) For $k=0,1$ the eigenvalues $
\th_{k1},\ldots,\th_{k,p+1}$ are distinct; and (b)~among the
values taken by $\rho_k[r_1,r_2;
\ell_{1},\ell_{2}]$ for $k,r_1,r_2=0,1$
and $1\leq\ell_{1},\ell_{2}\leq p$, the absolute value of $
\rho_k[r_1,r_2;\ell_{1},
\ell_{2}]$ equals 1 only when $r_1=r_2$ and $
\ell_{1}=\ell_{2}$.
\end{tabular}\hspace*{-34pt}
\end{equation}

Condition (\ref{eq315})(a) ensures that the eigenfunctions $\psi
_{k\ell}$ are well defined for $k=0,1$ and $\ell=1,\ldots,p$; and
(\ref{eq315})(b) guarantees that the quantities $\int_\cI (g-\mu
_{r_1}) \psi_{r_1\ell_{1}}$ and $\int_\cI (g-\mu_{r_2}) \psi
_{r_2\ell_{2}}$, which appear in the definition of $T_0(g)$ at (\ref
{eq216}), cannot be identical, except for a difference in means, unless
$r_1=r_2$ and $\ell_{1}=\ell_{2}$, thereby avoiding degeneracy.

The counterpart of Theorem \ref{Theorem1} for the quadratic
discriminant classifier is stated in the next theorem. Its proof is
given in the supplementary material [\citet{CarDelHal}].
%
\begin{theo}\label{Theorem3} Assume that (\ref{eq37})--(\ref
{eq39}) and (\ref{eq315}) hold. Then the error rate of the quadratic
discriminant, when the data in $\cD$ come from $\Pi_k$, admits the
expansion at (\ref{eq310}), but with different constants, where the
various terms have the properties stated immediately below that formula.
\end{theo}

Again the signs of $c_k$, $c_{k1}$ and $d_{k0} \nu_0^{-1}+d_{k1} \nu
_1^{-1}$, in (\ref{eq310}), are particularly diverse, and so the
conclusions reached in Section~\ref{Remark2} apply. Likewise, the
integration effect discussed in Section~\ref{Remark4} is also
observed. Here, as can be seen directly from (\ref{eq213}), the
estimator ${\widehat g}$ is integrated, and only the integral is
squared, not ${\widehat g}$ itself. The resulting integration effect
eliminates any term in $(mh)^{-1}$ from the analogue of the expansion
(\ref{eq310}) in this setting, although again this influence does not
carry over to the training data.

\section{Numerical illustrations}
\subsection{Simulated data}\label{secsimul}

To illustrate the impact of bandwidth on classification performance, we
generated data from several instances of model (\ref{eq21}), taking,
in each case, $m\kj=50$. Let $\phi_{\sigma}(x)$ denote the normal
density function with mean zero and standard deviation $\sigma$. We
considered the following cases, each with three different levels of
errors, which we refer to as noise versions 1, 2 and 3:

\begin{longlist}[(A):]
\item[(A):]
$g\kj(t)=\mu_k(t)+ (3 t+100)^{1/2} \{\cos(t/50)\}^k Z\kj$,
where
$\mu_0(t)=\phi_{10}(t-5)$,
$\mu_1(t)=\mu_0(t)+0.3\cos(t/5)+0.1$,
$Z\kj\sim U[-1/(30-10 k),1/(30-10 k)]$,
and
$\ep\kji\sim N(0,1/(4-2k)^2)$ (noise version 1),
$\ep\kji\sim N(0,2/(4-2k)^2)$ (noise version 2)
or $\ep\kji\sim N(0,4/(4-2k)^2)$ (noise version 3),
and $\pi_0=1/3$, $\pi_1=2/3$. Moreover, $X\kji=2i-1$, for
$i=1,\ldots,50$.
\end{longlist}

\begin{longlist}[(B):]
\item[(B):]
$g\kj(t)=\mu_k(t)+ (3 t+100)^{1/2} Z\kj$,
where
$\mu_0(t)
=30  \{0.2 \phi_4(t-5)+0.1 \phi_4(t-10)+0.4 \phi
_6(t-20)+0.4 \phi_6(t-35)+0.6 \phi_7(t-55)+0.6 \phi_7(t-80) \}$,
$\mu_1(t)=\mu_0(t)+4/\{(t-50)^2+10\},
$
$Z\kj\sim U[-1/(60+15 k),1/(60+15 k)]$,
$\ep\kji\sim\{\operatorname{Exp}(0.5)-2\}/(2+2k)$ (noise version 1),
$\ep\kji\sim\sqrt2 \{\operatorname{Exp}(0.5)-2\}/(2+2k)$ (noise version 2)
or $\ep\kji\sim\{\operatorname{Exp}(0.5)-2\}/(1+k)$ (noise version 3),
and $\pi_0=2/5$ and $\pi_1=3/5$. Moreover, $X\kji$ was as in (A).
\end{longlist}

\begin{longlist}[(C):]
\item[(C):]
$g_{0j}(t)=\mu_0(t)+ (3 t+100)^{1/2} Z_{0j}$,
$g_{1j}(t)=\mu_1(t)+ (t+5) Z_{1j}$,
where
$\mu_0(t)=15 \phi_{17}(t-65)\cos(t/7)$,
$\mu_1(t)=\mu_0(t)+5 \phi_{20}(t-50)$,
$Z\kj\sim U[-1/(50-10 k),1/(50-10 k)]$,
$\ep\kji\sim N(0,(4-k)^2/100)$ (noise version 1),
$\ep\kji\sim N(0,\allowbreak(4-k)^2/50)$ (noise version 2),
$\ep\kji\sim N(0,(4-k)^2/25)$ (noise version 3),
and $\pi_0=2/3$ and $\pi_1=1/3$.
Moreover, $X\kji$ was as in (A).
\end{longlist}

\begin{longlist}[(D)--(F):]
\item[(D)--(F):]
Same as (A) to (C) but with $X\kji=2i-1+T\kji$, where $T\kji\sim N(0,0.25)$.
\end{longlist}

We chose these examples to illustrate various features of the problems,
namely that the impact of smoothing may differ among classifiers,
and that in some cases, some classifiers perform better with more
smoothing and in other cases, they might perform better with less smoothing.

In each case, for $k=0,1$ and for several values of $n_{\mathrm{tr}}$,
we generated $100$ (resp.,~$n_{\mathrm{tr}}$) noisy test curves (resp.,
training curves) from model (\ref{eq21}), each of which came with
probability $\pi_k$ from $\Pi_k$. We constructed each classifier from
the training data, and applied it to the test data. To compute
${\widehat g}$ and ${\widehat g}\kj$, we compared three approaches for
selecting the bandwidths: no smoothing (NS), the standard plug-in (PI)
bandwidths $h_{\mathrm{PI}}$ and $h_{\mathrm{PI},kj}$ that estimate the
optimal bandwidth for estimation of the regression functions $g$ and
$g\kj$, which we computed using the \texttt{dpill} function in the R
package \texttt{KernSmooth}; see \citet{RupSheWan95}; and the
bandwidths $\gamma h_{\mathrm{PI}}$ and $\gamma_1 h_{\mathrm{PI},kj}$,
where $\gamma$ and $\gamma_1$ (and also the truncation parameter $p$ in
the case of the quadratic discriminant classifier) were chosen to
minimise the following cross-validation (CV) estimator of
classification error:
\[
\widehat\err=\frac{\widehat\pi_0}{n_0} \sum_{i=1}^{n_0}
I\{ \widehat\cC_{i0,-i}=1\}+ \frac{\widehat\pi_1}{n_1} \sum
_{i=1}^{n_1} I\{\widehat\cC_{i1,-i}=0\}
\]
with $\widehat\pi_0$ and $\widehat\pi_1$ denoting estimators of
$\pi_0$ and $\pi_1$ (we took $\widehat\pi_k=1/2$), and $\widehat
\cC_{ik,-i}$ being the estimator of the class label of the $i$th
training observation from group $k$, obtained from the classifier
constructed without using this observation.\vadjust{\goodbreak}

\begin{table}
\def\arraystretch{0.9}
\caption{Percentage of correctly classified
observations for the simulated data of Section
\protect\ref{secsimul}, using plug-in (PI) regression bandwidths,
bandwidths that minimise a crossvalidation (CV) estimate of
classification error, or without smoothing the noisy data (NS). The
three noise versions, in increasing order, are described in cases
\textup{(A)--(C)} in Section \protect\ref{secsimul}. Here ``Cent'' is the centroid
classifier (\protect\ref{eq29}), ``Cent sc.'' is the scaled centroid
classifier (\protect\ref{eq210}) and ``QDA'' is the quadratic discriminant
classifier (\protect\ref{eq213})} \label{tableshort}
\begin{tabular*}{\tablewidth}{@{\extracolsep{\fill}}lr ccc ccc ccc@{}}
\hline
&&\multicolumn{3}{c}{\textbf{Cent}}&\multicolumn{3}{c}{\textbf{Cent
sc.}}&\multicolumn{3}{c@{}}{\textbf{QDA}}\\[-4pt]
&& \multicolumn{3}{c}{\hrulefill} & \multicolumn{3}{c}{\hrulefill}
& \multicolumn{3}{c@{}}{\hrulefill}\\
&\multicolumn{1}{c}{$\bolds{n_{\mathrm{tr}}}$}&\textbf{CV} &\textbf{PI}
&\textbf{NS}&\textbf{CV} &\textbf{PI}&\textbf{NS}&\textbf{CV}
&\textbf{PI}&\multicolumn{1}{c@{}}{\textbf{NS}} \\
\hline
&&\multicolumn{9}{c@{}}{Case (A)}\\
[4pt]
Noise version 1&$50$& 82.9&74.1&84.0&91.8&73.2&92.0&95.1&94.1&53.5\\
Noise version 1&$100$& 84.4&74.9&84.8&92.6&74.1&92.6&97.6&94.8&67.6\\
[4pt]
Noise version 2&$50$& 77.7&69.6&78.1&94.3&70.4&94.4&91.0&89.3&49.1 \\
Noise version 2&$100$& 79.9&70.7&80.2&95.1&71.2&95.1&94.3&89.7&61.7\\
[4pt]
Noise version 3&$50$& 71.1&65.6&71.1&97.1&69.0&97.1&85.4&84.3&46.2 \\
Noise version 3&$100$& 73.7&66.8&74.1&97.9&69.4&97.9&89.9&84.1&58.4 \\
[4pt]
&&\multicolumn{9}{c@{}}{Case (B)}\\
[4pt]
Noise version 1&$50$
&63.2&60.1&65.7&96.3&78.7&96.5&77.1&74.3&65.8\\
Noise version 1&$100$
&65.5&61.5&66.8&96.8&80.0&96.8&81.8&76.3&73.0\\
[4pt]
Noise version 2&$50$
&61.5&58.6&64.6&96.3&80.6&96.4&76.9&74.1&65.2\\
Noise version 2&$100$
&62.6&58.7&64.4&96.7&81.3&96.7&81.3&75.0&72.4 \\
[4pt]
Noise version 3&$50$
&60.9&57.6&64.0&96.2&81.6&96.4&77.3&74.2&65.4\\
Noise version 3&$100$
&60.7&56.8&63.3&96.7&82.3&96.7&81.6&75.2&72.3 \\
[4pt]
&&\multicolumn{9}{c@{}}{Case (C)}\\
[4pt]
Noise version 1&$50$&61.5&60.8&60.8&88.7&89.2&87.4&84.8&83.7&82.0\\
Noise version 1&$100$&59.4&58.4&58.2&90.0&90.3&88.5&86.9&85.7&79.0\\
[4pt]
Noise version 2&$50$&61.3&60.2&60.3&87.3&87.9&82.8&81.9&81.2&82.4\\
Noise version 2&$100$&58.9&57.9&57.6&88.8&89.0&85.2&84.6&83.1&80.4\\
[4pt]
Noise version 3&$50$&61.0&59.7&59.3&85.2&85.4&71.2&80.5&79.9&79.6\\
Noise version 3&$100$&58.5&57.4&57.0&87.2&86.6&74.9&82.6&81.1&79.7\\
\hline
\end{tabular*}
\end{table}

For each configuration, we generated $B=200$ sets of training and test
samples. In Tables~\ref{tableshort} and \ref{tableshort2}, we
report the percentage of correctly classified test curves, averaged
over the $B$ replicates. Depending on the model, the classifier, and
the type of data (test or training), the cross-validation bandwidths
were either smaller or larger than the PI regression bandwidths,
illustrating the variety of settings already explained by our theory.
See Table B.1 
in Section B.3 in the supplementary material [\citet{CarDelHal}], where we report
the value of $\gamma$ and $\gamma_1$ averaged over the $B$
replicates. We can see from the table that in most cases, $\gamma$ was
smaller than $\gamma_1$, and both were usually smaller than $1$,
except in cases (C) and (F).

\begin{table}
\def\arraystretch{0.9}
\caption{Percentage of correctly classified
observations for the simulated data of Section \protect\ref
{secsimul}, using
plug-in (PI) regression bandwidths, bandwidths that minimise a
crossvalidation (CV) estimate of classification error, or without
smoothing the noisy data (NS). The three noise versions, in increasing
order, are described in cases \textup{(D)--(F)} in Section \protect\ref
{secsimul}. Here
``Cent'' is the centroid classifier (\protect\ref{eq29}), ``Cent
sc.'' is the
scaled centroid classifier (\protect\ref{eq210}) and ``QDA'' is the quadratic
discriminant classifier (\protect\ref{eq213})}
\label{tableshort2}
\begin{tabular*}{\tablewidth}{@{\extracolsep{\fill}}lr ccc ccc ccc@{}}
\hline
&&\multicolumn{3}{c}{\textbf{Cent}}&\multicolumn{3}{c}{\textbf{Cent
sc.}}&\multicolumn{3}{c@{}}{\textbf{QDA}}\\[-4pt]
&& \multicolumn{3}{c}{\hrulefill} & \multicolumn{3}{c}{\hrulefill}
& \multicolumn{3}{c@{}}{\hrulefill}\\
&\multicolumn{1}{c}{$\bolds{n_{\mathrm{tr}}}$}&\textbf{CV} &\textbf{PI}
&\textbf{NS}&\textbf{CV} &\textbf{PI}&\textbf{NS}&\textbf{CV}
&\textbf{PI}&\multicolumn{1}{c@{}}{\textbf{NS}} \\
\hline
&&\multicolumn{9}{c@{}}{Case (D)}\\
[4pt]
Noise version 1&$50$&
80.2&69.5&80.6&85.7&68.5&86.3&93.9&92.7&69.2\\
Noise version 1&$100$&
81.5&70.0&82.0&87.3&69.2&87.3&96.6&93.2&84.3\\
[4pt]
Noise version 2&$50$&
74.7&65.9&75.6&90.0&65.6&90.3&88.5&86.7&60.9\\
Noise version 2&$100$&
76.9&66.8&77.3&90.9&66.8&91.0&92.3&86.8&77.6\\
[4pt]
Noise version 3&$50$&
69.2&62.3&69.6&94.2&65.4&94.4&82.9&80.5&55.6\\
Noise version 3&$100$&
71.4&63.4&72.0&95.1&66.4&95.1&87.9&80.8&72.7\\
[4pt]
&&\multicolumn{9}{c@{}}{Case (E)}\\
[4pt]
Noise version 1&$50$&
65.0&61.9&67.3&94.8&79.0&95.0&77.0&73.1&71.2 \\
Noise version 1&$100$&
65.8&62.5&67.6&95.4&79.6&95.4&84.3&69.7&82.8\\
[4pt]
Noise version 2&$50$&
63.0&60.2&65.4&94.7&80.6&95.0&77.8&74.2&70.5 \\
Noise version 2&$100$&
63.4&59.9&64.9&95.4&81.4&95.5&84.8&69.8&82.4\\
[4pt]
Noise version 3&$50$&
61.3&59.2&64.3&94.6&81.4&94.9&77.9&74.4&69.8\\
Noise version 3&$100$&
61.8&58.2&63.1&95.4&82.5&95.5&84.5&71.5&81.8\\
[4pt]
&&\multicolumn{9}{c@{}}{Case (F)}\\
[4pt]
Noise version 1&$50$&
60.2&59.1&59.4&88.0&88.7&87.9&83.5&82.6&80.4\\
Noise version 1&$100$&
58.8&57.8&57.7&89.0&89.3&88.5&84.9&83.2&77.3\\
[4pt]
Noise version 2&$50$&
59.8&58.7&59.0&86.5&87.2&84.6&80.8&80.2&80.0\\
Noise version 2&$100$&
58.6&57.3&57.2&87.6&87.7&85.8&83.1&81.1&77.0\\
[4pt]
Noise version 3&$50$&
59.2&58.3&58.2&84.5&84.1&76.7&79.4&78.5&78.9\\
Noise version 3&$100$&
58.0&56.8&56.5&85.8&84.6&78.5&81.0&79.1&75.7\\
\hline
\end{tabular*}
\end{table}

As expected, we conclude from Tables~\ref{tableshort} and \ref
{tableshort2}, depending on the model and the classifier, the negative
impact of smoothing with the standard PI bandwidth can be quite
significant, indeed sometimes reducing the percentage of correctly
classified data by as much as 10\%. In cases (A) and (D), it is the
centroid classifier and its scaled version that are the most affected
by this inappropriate level of smoothing, whereas the quadratic
discriminant classifier is more robust against the level of smoothing.
In cases (B) and (E), the scaled centroid classifier and the quadratic
discriminant classifier are the most affected by inappropriate
smoothing. Cases (C) and (F) are more robust against smoothing; there,
all three versions (PI, CV and NS) of the data result in similar
classification performance, although overall the data smoothed by CV
result in slightly improved performance. Depending on the case, when
the noise level increases the impact of inappropriate bandwidth choice
can either increase or decrease.

\subsection{Real data}
We illustrate our findings on the ovarian cancer data set \texttt{8-7-02},
which concerns 253 patients (91 controls and 162 with\break ovarian cancer).
The data, which were produced to study the effect of robotic sample
handling, are
available from
\texttt{\href{http://home.ccr.cancer.gov/ncifdaproteomics/ppatterns.asp}{http://home.ccr.cancer.gov/}\break
\href{http://home.ccr.cancer.gov/ncifdaproteomics/ppatterns.asp}{ncifdaproteomics/ppatterns.asp}}.
In this example, the functions $X_i$ represent proteomic mass spectra
and $t\in[0,20\mbox{,}000]$ is the mass over charge ratio, $m/z$.
These raw curves are ideal for illustrating the negative impact that
systematically smoothing by standard methods can have, because in some
ranges of values of $t$, the spectra have considerable activity, and
the impact
of smoothing such data can be striking. We focus on one such ranges,
namely $t\in[200,500]$.

To assess the performance of classifiers on this data set, we randomly
and uniformly created $B=200$ pairs of (training sample, test sample),
where we took the training sample to be of size $n_{\mathrm{tr}}$ and
the test sample of size $253-n_{\mathrm{tr}}$, for $n_{\mathrm{tr}}=50$
and $n_{\mathrm{tr}}=100$. We also generated two more noise versions of
the data, adding to the $Y_{kji}$'s in both the test and training data,
noise $\varepsilon'_{kji} \sim\mathrm{N}(0,0.04)$ (noise version~1) or
$\varepsilon'_{kji}\sim\mathrm{N}(0,0.25)$ (noise version 2), where the
$\varepsilon'_{kji}$'s were totally independent.

For each version of the data (original data and noise versions 1 and
2), and for each pair of test and training sample, we constructed each
classifier from the training sample, and applied the classifier to the
test sample using either plug-in regression bandwidths to construct the
estimators ${\widehat g}$ and ${\widehat g}\kj$, or bandwidths
obtained by minimising the CV estimator of classification error defined
in Section~\ref{secsimul}, where we took $\widehat\pi_k=1/2$.

\begin{table}
\caption{Percentage of correctly classified
observations for the ovarian cancer data, using plug-in (PI) regression
bandwidths or bandwidths that minimise a crossvalidation (CV) estimate
of classification error. Here ``Cent'' is the centroid classifier
(\protect\ref{eq29}), ``Cent sc.'' is the scaled centroid classifier
(\protect\ref{eq210}) and ``QDA'' is the quadratic discriminant classifier
(\protect\ref{eq213})}
\label{tableovarian}
\begin{tabular*}{\tablewidth}{@{\extracolsep{\fill}}lrcccccc@{}}
\hline
&& \multicolumn{2}{c}{\textbf{Cent}} & \multicolumn{2}{c}{\textbf{Cent sc.}} &
\multicolumn{2}{c@{}}{\textbf{QDA}} \\[-4pt]
&& \multicolumn{2}{c}{\hrulefill} & \multicolumn{2}{c}{\hrulefill} &
\multicolumn{2}{c@{}}{\hrulefill} \\
\textbf{Data}&\multicolumn{1}{c}{$\bolds{n_{\mathrm{tr}}}$}& \textbf{CV}
&\textbf{PI} & \textbf{CV}&\textbf{PI} & \textbf{CV}&\textbf{PI} \\
\hline
Original data& $50$& 90.60 & 80.25 & 90.05 & 78.79 & 93.32 & 89.69\\
Original data&$100$&90.43 & 80.96 & 90.00 & 79.96 & 98.58 & 98.86 \\
[4pt]
Noisy version 1& $50$& 88.07 & 75.19 & 87.83 & 74.23 & 78.03 & 68.50\\
Noisy version 1& $100$&87.58 & 76.76 & 88.54 & 76.27 & 91.48 & 90.97 \\
[4pt]
Noisy version 2& $50$& 76.15 & 66.57 & 76.65 & 66.09 & 56.91 & 48.54\\
Noisy version 2&$100$&81.97 & 67.55 & 81.91 & 67.64 & 77.62 & 66.49 \\
\hline
\end{tabular*}
\end{table}

Table~\ref{tableovarian} reports the percentage, averaged over the
$B$ pairs of samples, of correctly classified observations from the
test samples.
The table indicates very clearly that smoothing the data using the
plug-in regression bandwidths degraded the quality of the two versions
of the centroid classifier by about 10\%, and a similar phenomenon was
observed for the quadratic discriminant classifier when the training
sample was small and when the data were noisy.

\begin{appendix}\label{app}
\section*{Appendix: Proof of Theorem 1}\label{app1}

\subsection{Preliminary results}\label{app11}
Define
\begin{eqnarray*}
\De_\ell(x)&=&{1\over mh} \sum
_{i=1}^m \ep_i \biggl(
{x-X_i\over h} \biggr)^{ \ell} K \biggl({x-X_i\over h}
\biggr),
\\
W_\ell(x)&=&{1\over mh} \sum_{i=1}^m
\biggl[\int_x^{X_i}\bigl\{g''(t)-g''(x)
\bigr\} (X_i-t) \,dt \biggr] \biggl({x-X_i\over h}
\biggr)^{ \ell} K \biggl({x-X_i\over h} \biggr).\vadjust{\goodbreak}
\end{eqnarray*}
With $U_\ell$ and $V_\ell$ given by
(\ref{eq24}) and
(\ref{eq26}), and using the model at (\ref{eq22}) and the exact
form of the remainder in Taylor's theorem, we can write:
\begin{eqnarray*}
V_\ell(x)&=&{1\over mh} \sum_{i=1}^m
\bigl\{g(X_i)+\ep_i\bigr\} \biggl({x-X_i\over h}
\biggr)^{ \ell} K \biggl({x-X_i\over h} \biggr)
\\
&=&{1\over mh} \sum_{i=1}^m
\biggl[g(x)+(X_i-x) g'(x)+{ {1\over2}}(X_i-x)^2
g''(x)+\ep_i \biggr]
\\
&&\hspace*{30.5pt}{}\times \biggl({x-X_i\over h} \biggr)^{ \ell} K
\biggl({x-X_i\over h} \biggr) +W_\ell(x)
\\
&=&g(x) U_\ell(x)-h g'(x) U_{\ell+1}(x) +{
{1\over2}}h^2 g''(x)
U_{\ell+2}(x)+\De_\ell(x)+W_\ell(x).
\end{eqnarray*}
Assuming, without loss of generality, that $K$ is supported on $[-1,1]$,
\[
\bigl|W_\ell(x)\bigr| \leq h^2 \Bigl\{\sup_{t\in\cI: |t-x|\leq h}
\bigl|g''(t)-g''(x)\bigr| \Bigr\}
{1\over mh} \sum_{i=1}^m K
\biggl({x-X_i\over h} \biggr) \leq h^2 U_0(x) Q,
\]
where $Q=\sup_{s,t\in\cI: |s-t|\leq h} |g''(s)-g''(t)|$. Now,
\[
{\widehat g}={U_2 V_0-U_1 V_1\over U_2 U_0-U_1^2} =g+{ {1\over2}}h^2
g'' {U_2^2-U_1 U_3\over U_2 U_0-U_1^2}+\De +{U_2 W_0-U_1 W_1\over U_2 U_0-U_1^2},
\]
where
$
\De= (U_2 \De_0-U_1 \De_1 ) / (U_2 U_0-U_1^2
).
$
Therefore, since $|U_\ell|\leq U_0$ for each $\ell\geq0$,
%
%
\begin{equation}\label{eqP1}
\biggl|{\widehat g}- \biggl(g+{ {1\over2}}h^2
g'' {U_2^2-U_1 U_3\over U_2
U_0-U_1^2} +\De \biggr) \biggr| \leq
{2 Q h^2 U_0^2 \over U_2 U_0-U_1^2},
\end{equation}
uniformly on $\cI$.

Similarly, defining $Q\kj=\sup_{s,t\in\cI: |s-t|\leq h} |g\kj
''(s)-g\kj''(t)|$, and
\begin{eqnarray*}
\De_{kj\ell}(x)&=&{1\over m\kj h_1} \sum
_{i=1}^{m\kj} \ep\kji \biggl({x-X\kji\over h_1}
\biggr)^{ \ell} K \biggl({x-X\kji\over h_1} \biggr),
\\
\De_{kj}&=&{U_{kj2} \De_{kj0}-U_{kj1} \De_{kj1}\over
U_{kj2} U_{kj0}-U_{kj1}^2},
\end{eqnarray*}
where $U_{kj\ell}$ is as at (\ref{eq25}), we have, uniformly on $\cI$
%
%
\begin{equation}\label{eqP2}\qquad
\biggl|{\widehat g}\kj- \biggl(g\kj+{ {1\over2}}h_1^2
g\kj'' {U_{kj2}^2-U_{kj1} U_{kj3}\over
U_{kj2} U_{kj0}-U_{kj1}^2} +\De\kj \biggr) \biggr| \leq
{2 Q\kj h_1^2 U_{kj0}^2 \over U_{kj2} U_{kj0}-U_{kj1}^2}.
\end{equation}

Define
%
%
\begin{equation}\label{eqP3}
\bDe_k={1\over n_k} \sum_{j=1}^{n_k}
\De\kj
\end{equation}
and recall that $\ka_2=\int u^2 K(u) \,du$. We shall derive the
following result in Section~\ref{app16}:

\begin{Lem}\label{LemmaP1} Under the conditions of Theorem \ref
{Theorem1}, for some $C_1>0$, all $C_2>0$ and $k=0,1$,
\begin{eqnarray*}
P_k \biggl(\sup_\cI \biggl|{U_2^2-U_1 U_3\over U_2 U_0-U_1^2}-
\ka _2 \biggr| >n^{-C_1} \biggr)&=&O \bigl(n^{-C_2} \bigr),
\\
P_k \biggl(\max_{j=1,\ldots,n_k} \sup_\cI
\biggl|{U_{kj2}^2-U_{kj1} U_{kj3}\over
U_{kj2} U_{kj0}-U_{kj1}^2}-\ka_2 \biggr| >n^{-C_1} \biggr)&=&O
\bigl(n^{-C_2} \bigr)
\end{eqnarray*}
as $n\rai$, and for some $C_3>0$, all $C_2>0$ and $k=0,1$,
\begin{eqnarray*}
P_k \biggl(\sup_\cI{U_0^2 \over U_2 U_0-U_1^2}>C_3
\biggr) &=& O \bigl(n^{-C_2} \bigr),
\\
P_k \biggl(\max_{j=1,\ldots,n_k} \sup_\cI
{U_{kj0}^2 \over U_{kj2} U_{kj0}-U_{kj1}^2}>C_3 \biggr) &=&O \bigl(n^{-C_2}
\bigr).
\end{eqnarray*}
Furthermore, defining $M\ssum=\min_{k=1,2}(\sum_j m\kj)$, we have for
all $C_2,C_4>0$,
%
%
\begin{eqnarray}\label{eqP4}
&&P_k \Bigl\{\sup_\cI|\De|>n^{C_4}
(mh)^{-1/2} \Bigr\}
\nonumber
\\
&&\quad{}+\max_{k=0,1} P_k \Bigl\{\sup_\cI|
\bDe_k|>n^{C_4} (M\ssum h)^{-1/2} \Bigr\} \\
&&\qquad=O
\bigl(n^{-C_2} \bigr).\nonumber
\end{eqnarray}
\end{Lem}

\subsection{\texorpdfstring{Initial calculation of $\err_k$}
{Initial calculation of err k}}\label{app12}
Let $\cG_1$ denote the sigma-field generated by the random variables
introduced in Section~\ref{sec2}, and the random functions $g\kj$,
but excluding $g$. Specifically, $\cG_1$ is the sigma-field generated
by $g\kj$, $X\kji$ and $\ep\kji$ for $1\leq i\leq m\kj$, $1\leq
j\leq n_k$ and $k=0,1$, and by $X_i$ and $\ep_i$ for $1\leq i\leq m$.
Recall that $\bow_k$ is $>$ or $\leq$ according as $k=0$ or $k=1$,
respectively, and recall formula (\ref{eq313}) for the
statistic $S({\widehat g})$.

Under the assumption that the new data set $\cD$ comes from $\Pi_k$,
and conditional on $\cG_1$, ${\widehat g}$ is a Gaussian\vspace*{1pt} process with
mean $\hal_k=E_k({\widehat g}\mi\cG_1)$ and covariance
function~$\hGa_k$, say. In this notation,
%
%
\begin{equation}\label{eqP5}
\err_k\equiv E_k\bigl[P_k\bigl\{S({\widehat
g})\bow_k0\mi\cG_1\bigr\}\bigr] =E_k\bigl\{
\Psi_k(-\hbe_k/\hsi_k)\bigr\},
\end{equation}
where, by (\ref{eq313}),
%
%
\begin{eqnarray}\qquad
\label{eqP6}\hbe_k&=&E_k\bigl\{S({\widehat g})\mi\cG_1
\bigr\} =\int_\cI (\hmu_1-\hmu_0)
\bigl\{2 \hal_k-(\hmu_0+\hmu_1)\bigr\},
\\
\label{eqP7}
\hsi_k^2&=&\var\bigl\{S({\widehat g})\mi
\cG_1\bigr\} \nonumber\\
&=&4\int_\cI \int
_\cI \bigl\{\hmu_1(x_1)-
\hmu_0(x_1)\bigr\} \bigl\{\hmu_1(x_2)-
\hmu_0(x_2)\bigr\}
\\
&&\hspace*{29pt}{}\times \hGa_k(x_1,x_2) \,dx_1
\,dx_2.\nonumber
\end{eqnarray}
The probability on the left-hand side of (\ref{eqP5}) equals the
chance that, when $\cD$ comes from $\Pi_k$, the classifier based on
$S({\widehat g})$ makes an error and assigns $\cD$ to the other
population.

\subsection{\texorpdfstring{Approximations to $\hal_k$, $\hbe_k$ and $\hsi_k$}
{Approximations to alpha k, beta k and sigma k}}\label{app13}
In view of (\ref{eqP1}),
%
%
\begin{equation}\label{eqP8}
\biggl|\hal_k- \biggl(\mu_k+{ {1\over2}}h^2
\mu_k'' {U_2^2-U_1 U_3\over U_2 U_0-U_1^2} +\De
\biggr) \biggr| \leq{2 E_k(Q) h^2 U_0^2 \over U_2 U_0-U_1^2}.
\end{equation}
Noting that, for random variables $A_1$, $A_2$, $B_1$ and $B_2$,
$
|{\cov}(A_1+A_2,B_1+B_2)-\cov(A_1,A_2) |
\leq|{\cov}(B_1,B_2)|+|{\cov}(A_1,B_2)|
+|{\cov}(B_1,A_2)|
$
where the covariances are interpreted conditionally on $\cG_1$, we
deduce from (\ref{eqP1}) that for a constant $C_4>0$,
%
%
\begin{eqnarray}\label{eqP9}\quad
&&\sup_{x_1,x_2\in\cI} \biggl|\hGa_k(x_1,x_2)-
\biggl\{G_k(x_1,x_2) +{
{1\over2}}h^2 G_k^{(0,2)}(x_1,x_2)
{U_2^2-U_1 U_3\over U_2 U_0-U_1^2} (x_2)
\nonumber
\\
&&\qquad\hspace*{123pt}{}+{ {1\over2}}h^2
G_k^{(2,0)}(x_1,x_2)
{U_2^2-U_1 U_3\over U_2 U_0-U_1^2} (x_1) \biggr\} \biggr|
\\
&&\qquad\leq C_4 h^2 \bigl\{h^2+E_k
\bigl(Q+Q^2 \bigr) \bigr\} \sup_\cI \biggl(1+
{U_0^2\over U_2 U_0-U_1^2} \biggr)^{ 2},\nonumber
\end{eqnarray}
where we define $G_k^{(j_1,j_2)}(x_1,x_2)=\partial
^{j_1+j_2}G_k(x_1,x_2)/\partial x_1^{j_1} \,\partial x_2^{j_2}$. (Recall
that $G_k$ denotes the covariance of the Gaussian process $g$ when the
data $\cD$ are drawn from~$\Pi_k$.)

With $\bg_k$ defined as at (\ref{eq31}), and defining $\bDe_k$ as
at (\ref{eqP3}), we have, in view of~(\ref{eqP2}), Lemma \ref
{LemmaP1} and (\ref{eq38})(b), the result
%
%
\begin{equation}\label{eqP10}
P_k \biggl\{\sup_\cI \biggl|\hmu_k -
\biggl(\bg_k+{ {1\over2}}h_1^2
\ka_2 \bg_k'' +
\bDe_k \biggr) \biggr|>n^{-C_1} h_1^2
\biggr\} =O \bigl(n^{-C_2} \bigr)
\end{equation}
for some $C_1>0$ and all $C_2>0$. Using Rosenthal's inequality, it can
be proved from (\ref{eq38}) and (\ref{eq39})(c) that, for some
$C_1>0$ and all $C_2>0$,
%
%
\begin{equation}\label{eqP11}
P_k \Bigl(\sup_\cI \bigl|\bg_k''-
\mu_k'' \bigr|>n^{-C_1} \Bigr) =O
\bigl(n^{-C_2} \bigr).
\end{equation}
Together, (\ref{eqP10}) and (\ref{eqP11}) imply that
%
%
\begin{equation}\label{eqP12}\quad
P_k \biggl\{\sup_\cI \biggl|\hmu_k -
\biggl(\bg_k+{ {1\over2}}h_1^2
\ka_2 \mu_k'' +
\bDe_k \biggr) \biggr|>n^{-C_1} h_1^2
\biggr\} =O \bigl(n^{-C_2} \bigr).
\end{equation}

Define $H^2=h^2+h_1^2$,
%
%
\begin{eqnarray}\label{eqP13}
\be_k&=&\int_\cI  \biggl\{\bg_1-
\bg_0 +{ {1\over2}}h_1^2
\ka_2 \bigl(\mu_1''-
\mu_0''\bigr) +\bDe_1-
\bDe_0 \biggr\} \nonumber\\
&&\hspace*{8.5pt}{}\times\biggl\{2 \mu_k-(\bg_0+
\bg_1)+h^2 \ka_2 \mu_k''
\nonumber
\\
&&\hspace*{25pt}{}
-{ {1\over2}}h_1^2 \ka_2
\bigl(\mu_0''+\mu_1''
\bigr) +2 \De-(\bDe_0+\bDe_1) \biggr\},
\\
\tsi_k^2&=&4\int_\cI \int
_\cI  \biggl\{\bg_1-\bg_0+{
{1\over2}}h_1^2 \ka_2 \bigl(
\mu_1''-\mu_0''
\bigr) +\bDe_0-\bDe_1 \biggr\}(x_1)
\nonumber
\\
&&\hspace*{29pt}{} \times \biggl\{\bg_1-\bg_0+{ {1\over2}}h_1^2
\ka_2 \bigl(\mu_1''-
\mu_0''\bigr) +\bDe_0-
\bDe_1 \biggr\}(x_2)
\nonumber
\\
&&\hspace*{29pt}{} \times \biggl[G_k(x_1,x_2) +{
{1\over2}}h^2 \ka_2 \bigl
\{G_k^{(2,0)}(x_1,x_2)
+G_k^{(0,2)}(x_1,x_2) \bigr\} \biggr].
\nonumber
\end{eqnarray}
Combining Lemma \ref{LemmaP1}, (\ref{eqP5})--(\ref{eqP9}) and
(\ref{eqP12}), we deduce that, for some $C_1>0$ and all $C_2>0$,
%
%
\begin{eqnarray}\label{eqP14}
P_k \bigl( |\hbe_k-\be_k |>n^{-C_1}
H^2 \bigr) &=&O \bigl(n^{-C_2} \bigr),\nonumber\\[-8pt]\\[-8pt]
P_k \bigl(
\bigl|\hsi_k^2-\tsi_k^2
\bigr|>n^{-C_1} H^2 \bigr) &=&O \bigl(n^{-C_2} \bigr).\nonumber
\end{eqnarray}
Observe from (\ref{eqP13}) that $\be_k=\be_{k0}+b_{k1}+\be
_{k1}+\be_{k2}+\bDe_2$, where $\be_{k0}$ is as at~(\ref{eq33}),
%
%
\begin{eqnarray}\label{eq34}
b_{k1} &=& \ka_2 \int_\cI (
\mu_1-\mu_0) \bigl(h^2\mu_k''-
h_1^2 \mu_{1-k}''
\bigr),
\\
\be_{k1}&=&\int_\cI (\bg_1-
\bg_0) \bigl\{2 \De-(\bDe_0+\bDe_1)\bigr\}\nonumber\\
&&{} +
\int_\cI \bigl\{2 \mu_k-(\bg_0+
\bg_1)\bigr\} (\bDe_1-\bDe_0),
\nonumber
\\
\be_{k2}&=&\int_\cI \bigl\{2 \De-(
\bDe_0+\bDe_1)\bigr\} (\bDe_1-
\bDe_0)
\nonumber
\end{eqnarray}
and $\bDe_2=\be_k-(\be_{k0}+b_{k1}+\be_{k1}+\be_{k2})$. Using
(\ref{eqP4}) it can be shown that, for some $C_1>0$ and all $C_2>0$,
and when $\ell=2$,
%
%
\begin{equation}\label{eqP15}
P_k \bigl(|\bDe_\ell|>n^{-C_1} H^2
\bigr) =O \bigl(n^{-C_2} \bigr).
\end{equation}
Hence, noting the first result in (\ref{eqP14}), we have:
%
%
\begin{equation}\label{eqP16}
P_k \bigl\{ \bigl|\hbe_k-(\be_{k0}+b_{k1}+
\be_{k1}+\be_{k2}) \bigr|>n^{-C_1} H^2 \bigr\}
=O \bigl(n^{-C_2} \bigr).
\end{equation}
Recall the definitions of $\si_k^2$ and $\tau_k^2$ at (\ref{eq35})
and (\ref{eq36}), and put
%
%
\begin{eqnarray}\label{eqP17}\qquad
\si_{k0}&=&2 h^2 \ka_2\int_\cI
\int_\cI (\bg_1-\bg_0)
(x_1) (\bg_1-\bg_0) (x_2)
\nonumber\\[-8pt]\\[-8pt]
&&\hspace*{48.7pt}{} \times \bigl\{G_k^{(2,0)}(x_1,x_2)
+G_k^{(0,2)}(x_1,x_2) \bigr\}
\,dx_1 \,dx_2,\nonumber
\\
\label{eqP18}
\si_{k1}&=&4 h_1^2 \ka_2\int
_\cI \int_\cI (\bg_1-
\bg_0) (x_1) (\mu_1-\mu_0)''(x_2)
G_k(x_1,x_2) \,dx_1
\,dx_2
\end{eqnarray}
and $\bDe_3=\tsi_k^2-(\si_k^2+\si_{k0}+\si_{k1})$. Thus, $\bDe_3$
is the term in $\bDe_0$ and $\bDe_1$ that arises when $\tsi_k^2$ is
expanded. Using (\ref{eqP4}) it can be proved that (\ref{eqP15})
holds when $\ell=3$. Moreover, $\hsi_k^2$ can be written as
%
%
\begin{equation}\label{eqP19}
\hsi_k^2 =\si_k^2+
\si_{k0}+\si_{k1}+\bDe_3+\bDe_4,
\end{equation}
where, in view of the second part of (\ref{eqP14}), (\ref{eqP15})
holds in the case $\ell=4$ and for some $C_1>0$ and all $C_2>0$.

Define $\tau_{k\ell}$ to be equal to $\si_{k\ell}$, at (\ref
{eqP17}) and (\ref{eqP18}), when $\bg_0$ and $\bg_1$ on the
respective right-hand sides are replaced by $\mu_0$ and $\mu_1$. Then
for $k=0,1$ and $\ell=0,1$, noting property (\ref{eq39})(c) on the
rates of increase of $n_0$ and $n_1$, it can be shown that for some $C_1>0$,
%
%
\begin{equation}\label{eqP20}
P_k \bigl(|\si_{k\ell}-\tau_{k\ell}|>n^{-C_1}
h_\ell^2 \bigr) =O \bigl(n^{-C_2} \bigr)
\end{equation}
for all $C_2>0$, where we define $h_0=h$. Therefore, if $C_1>0$ is
sufficiently small,
%
%
\begin{equation}\label{eqP21}
\max_{k=0,1} \max_{\ell=0,1} P_k
\bigl(|\si_{k\ell}|>n^{-C_1} \bigr) =O \bigl(n^{-C_2}
\bigr)
\end{equation}
for all $C_2>0$.

\subsection{\texorpdfstring{Approximation to $\hsi_k^{-1}$}{Approximation to sigma k -1}}\label{app14}

In the notation at (\ref{eqP19}),
\[
{1\over\hsi_k} ={1\over\tau_k} \biggl(1+
{\si_k^2-\tau_k^2\over\tau_k^2} +{\si_{k0}+\si_{k1}+\bDe_3+\bDe_4\over\tau_k^2} \biggr)^{ -1/2}
=s_k(\infty),
\]
where, for $0\leq r\leq\infty$,
\[
s_k(r)={1\over\tau_k} \sum_{j=0}^r
\sum_{\ell=0}^j \pmatrix{-{ \dfrac{1}{2}}
\vspace*{2pt}\cr j} \pmatrix{j\cr\ell} \biggl({\si_k^2-\tau_k^2\over\tau_k^2} \biggr)^{j-\ell}
\biggl({\si_{k0}+\si_{k1}+\bDe_3+\bDe_4\over\tau_k^2} \biggr)^{
\ell}.
\]
We claim that the infinite series defined by $s_k(\infty)$ converges
with probability $1-O(n^{-C_2})$ for all $C_2>0$. To appreciate why,
note that, by (\ref{eq38}) and (\ref{eq39})(c), there exists
$C_1>0$ such that
\[
P_k \bigl( \bigl|\si_k^2-\tau_k^2
\bigr|>n^{-C_1} \bigr)=O \bigl(n^{-C_2} \bigr)
\]
for all $C_2>0$. Combining this property, (\ref{eqP15}) for $\ell=3$
and 4, and (\ref{eqP21}), we deduce that, for some $C_1>0$ and all $C_2>0$,
\[
P_k \biggl( \biggl|{\si_k^2-\tau_k^2\over\tau_k^2} \biggr| + \biggl|{\si_{k0}+\si_{k1}+\bDe_3+\bDe_4\over\tau_k^2}
\biggr|\leq n^{-C_1} \biggr) =1-O \bigl(n^{-C_2} \bigr).
\]
Therefore, if $C_3>0$ is given then $r_0=r_0(C_3)\geq1$ can be chosen
so large that, whenever $r_0\leq r\leq\infty$,
$
P_k \{ |\hsi_k^{-1}-s_k(r) |>n^{-C_3} \}=O
(n^{-C_2} )
$
for all $C_2>0$. Using this property and (\ref{eqP15}), again for
$\ell=3$ and 4; and employing too (\ref{eqP20}); we see that for
some $C_1>0$ and all $C_2>0$, if $r_0$ is chosen sufficiently large,
%
%
\begin{equation}\label{eqP22}
P_k \bigl\{ \bigl|\hsi_k^{-1}-t_k(r)
\bigr|>n^{-C_1} H^2 \bigr\} =O \bigl(n^{-C_2} \bigr)
\end{equation}
for $r\geq r_0$, where
%
%
\begin{equation}\label{eqP23}\qquad
t_k(r)={1\over\tau_k} \sum_{j=0}^r
\sum_{\ell=0}^{\min(j,1)} \pmatrix{-{ \dfrac{1}{2}}
\vspace*{2pt}\cr j} \pmatrix{j\cr\ell} \biggl({\si_k^2-\tau_k^2\over\tau_k^2} \biggr)^{j-\ell}
\biggl({\tau_{k0}+\tau_{k1}\over\tau_k^2} \biggr)^{ \ell}.
\end{equation}

\subsection{\texorpdfstring{Approximation to $E_k\{\Psi_k(-\hbe_k/\hsi_k)\}$}
{Approximation to E k\{Psi k(-beta k/sigma k)\}}}\label{app15}

Let $C_1>0$ and let $\ell_0\geq0$ be an integer. With $U_{kj\ell}$
defined as at (\ref{eq25}), let $\cE$ denote the event
\[
\cE=\cE(C_1,\ell_0) = \Bigl\{\max_{1\leq\ell\leq\ell_0}
\max_{j=1,\ldots,n_k} \sup_{x\in\cI} \bigl|U_{kj\ell}(x)-
\ka_\ell f_X(x) \bigr|\leq n^{-C_1} \Bigr\},
\]
where $\ka_\ell=\int u^\ell K(u) \,du$ and hence vanishes for
odd $\ell$, since by (\ref{eq39})(b), $K$ is symmetric. It will be
proved in Section~\ref{app16} that, for some $C_1>0$ and each $\ell
_0\geq0$,
%
%
\begin{equation}\label{eqP24}
P_k\bigl\{\cE(C_1,\ell_0)\bigr\}=1-O
\bigl(n^{-C_2} \bigr) \qquad\mbox{for all } C_2>0.
\end{equation}
If $\cE(C_1,\ell_0)$ holds for an $\ell_0\geq2$ then, if
$0<C_1'<C_1$, there exists a nonrandom integer $n_0\geq1$ such that
the event $\cE_1=\cE_1(C_1')$, defined by
%
%
\begin{equation}\label{eqP25}\qquad\quad
\cE_1= \Bigl\{\max_{j=1,\ldots,n_k} \sup
_{x\in\cI} \bigl|U_{kj2}(x) U_{kj0}(x)-U_{kj1}(x)^2-
\ka_2 f_X(x)^2 \bigr| \leq n^{-C_1'} \Bigr
\}
\end{equation}
holds for all $n\geq n_0$.

Let $I=I(\cE)$ denote the indicator of $\cE$. In view of (\ref{eqP24}),
%
%
\begin{equation}\label{eqP26}
E_k\bigl\{\Psi_k(-\hbe_k/\hsi_k)
\bigr\} =E_k\bigl\{I \Psi_k(-\hbe_k/
\hsi_k)\bigr\}+O \bigl(n^{-C_2} \bigr)
\end{equation}
for all $C_2>0$, and so to approximate the term on the left-hand side
of (\ref{eqP26}) we may develop an approximation to the first term on
the right-hand side.

Let $\cG_2$ denote the sigma-field generated by the random variables
$X_i$ for \mbox{$1\leq i\leq m$}, and by $X\kji$ and the functions $g\kji$
for $1\leq i\leq m\kj$, $1\leq j\leq n_k$ and $k=0,1$ (i.e.,
generated by everything except $g$ and the experimental errors $\ep_i$
and~$\ep\kji$). The quantities $I$, $t_k(r)$ at (\ref{eqP23}), $\be
_{k0}$ at (\ref{eq33}), and $b_{k1}$ at (\ref{eq34}) are all $\cG
_2$-measurable. Therefore, using (\ref{eqP16}) and (\ref{eqP22}),
and noting that $\Psi_k$ is an analytic function with all derivatives
uniformly bounded, we obtain
%
%
\begin{eqnarray} \label{eqP27}
&&E_k\bigl\{I \Psi_k(-\hbe_k/
\hsi_k)\bigr\}
\nonumber
\\
&&\qquad=E_k \bigl(E_k \bigl[I \Psi_k \bigl\{-(
\be_{k0}+b_{k1}+\be_{k1}+\be _{k2})
t_k(r) \bigr\}\mid\cG_2 \bigr] \bigr)+o
\bigl(H^2 \bigr)
\nonumber
\\
&&\qquad=E_k\bigl[I \Psi_k\bigl\{-\be_{k0}
t_k(r)\bigr\}\bigr] -b_{k1} \tau_k^{-1}E_k
\bigl[I \Psi_k'\bigl\{-\be_{k0}
t_k(r)\bigr\}\bigr]
\\
&&\qquad\quad{} -\tau_k^{-1}E_k \bigl[E_k(
\be_{k2}\mi\cG_2) I \Psi_k'\bigl\{-
\be_{k0} t_k(r)\bigr\} \bigr]
\nonumber
\\
&&\qquad\quad{} +{ \tfrac{1}{2}}\tau_k^{-2}E_k
\bigl[E_k \bigl(\be_{k1}^2 \mid \cG _2
\bigr) I \Psi_k''\bigl\{-\be_{k0}
t_k(r)\bigr\} \bigr]
\nonumber
\\
&&\qquad\quad{} +O \bigl\{(mh)^{-2}+(M\ssum h_1)^{-2} \bigr\}
+o \bigl(H^2 \bigr).
\end{eqnarray}
Here we have used the properties $E_k(\be_{k1}\mi\cG_2)=0$,
$E_k|t_k(r)-\tau_k^{-1}|=O(n^{-C})$ for some $C>0$,
\[
E_k \bigl[E_k \bigl(\be_{k\ell_1}^{\ell_2} \mid
\cG_2 \bigr) I \Psi_k^{(\ell_2)}\bigl\{-
\be_{k0} t_k(r)\bigr\} \bigr] =O \bigl\{(mh)^{-2}+(M
\ssum h_1)^{-2} \bigr\}
\]
for $\ell_2\geq3$ if $\ell_1=1$, and for $\ell_2\geq2$ if $\ell
_1=2$, and
\[
\bigl|E_k \bigl[E_k (\be_{k1} \be_{k2} \mid
\cG_2 ) I \Psi_k''\bigl\{-
\be_{k0} t_k(r)\bigr\} \bigr] \bigr| =O \bigl
\{(mh)^{-2}+(M\ssum h_1)^{-2} \bigr\}.
\]
Further, we have used the fact that the event $\cE_1$, defined at
(\ref{eqP25}), obtains whenever $I\neq0$.

In addition,
%
%
\begin{eqnarray}\label{P28}\qquad
{ {1\over4}}E_k\bigl[E_k \bigl(
\be_{k1}^2 \mid \cG_2 \bigr) I \bigr]
&=&E_k \biggl\{I\int_\cI (\bg_0-
\bg_1) \De \biggr\}^{ 2}
\nonumber
\\
&&{}+E_k \biggl\{I\int_\cI (\bg_0-
\mu_k) \bDe_0 \biggr\}^{ 2} +E_k
\biggl\{I\int_\cI (\bg_1-\mu_k)
\bDe_1 \biggr\}^{ 2}
\\
&=&O \bigl(m^{-1} \bigr),\nonumber
\end{eqnarray}
that
%
%
\begin{eqnarray}\label{eqP29}
&&E_k \bigl[E_k(\be_{k2}\mi\cG_2) I
\Psi_k'\bigl\{-\be_{k0} t_k(r)
\bigr\} \bigr]
\nonumber
\\
&&\qquad=(-1)^{k+1} \phi(b_{k0}/\tau_k) \int
_\cI E_k \bigl[I \bigl\{E_k \bigl(
\bDe_0^2 \mid \cG_2 \bigr) -E_k
\bigl(\bDe_1^2 \mid \cG_2 \bigr) \bigr\}
\bigr]\nonumber\\[-8pt]\\[-8pt]
&&\qquad\quad{}+o \bigl\{(\nu _0h_1)^{-1} \bigr\}
\nonumber\\
&&\qquad={\ka\over h_1} \bigl( \si_{\ep0}^2
\nu_0^{-1}-\si_{\ep1}^2
\nu_1^{-1} \bigr) (-1)^{k+1} \phi(b_{k0}/
\tau_k) \int_\cI f_X^{-1}
+o \bigl\{(\nu_0h_1)^{-1} \bigr\}\nonumber
\end{eqnarray}
and that
%
%
\begin{equation}\label{eqP30}\qquad
b_{k1} \tau_k^{-1}E_k\bigl[I
\Psi_k'\bigl\{-\be_{k0} t_k(r)\bigr
\}\bigr] =b_{k1} \tau_k^{-1}(-1)^{k+1}
\phi(b_{k0}/\tau_k) +o \bigl(H^2 \bigr),
\end{equation}
where $b_{k0}$ and $b_{k1}$ are as at (\ref{eq33}) and (\ref
{eq34}), $\phi$ is the standard normal density, and we have used the
fact that $\Psi_k'=(-1)^{k+1} \phi$. Combining (\ref{eqP24}) and
(\ref{eqP26})--(\ref{eqP30}), and taking $r$ sufficiently large
(but fixed), we deduce that
%
%
\begin{eqnarray}\label{eqP31}\qquad
E_k\bigl\{\Psi_k(-\hbe_k/\hsi_k)
\bigr\} &=&E_k\bigl[\Psi_k\{-\be_{k0}/
\si_k\}\bigr] -b_{k1} \tau_k^{-1}(-1)^{k+1}
\phi(b_{k0}/\tau_k)
\nonumber
\\
&&{} -{\ka\over\tau_kh_1} \bigl(\si_{\ep0}^2
\nu_0^{-1}-\si_{\ep1}^2
\nu_1^{-1} \bigr) (-1)^{k+1} \phi(b_{k0}/
\tau_k)\int_\cI f_X^{-1}
\\
&&{} +O \bigl\{m^{-1}+(mh)^{-2} \bigr\} +o \bigl
\{H^2+(\nu_0h_1)^{-1} \bigr\}.\nonumber
\end{eqnarray}
Result (\ref{eq310}) follows from (\ref{eqP5}) and (\ref{eqP31}).

\subsection{\texorpdfstring{Proof of Lemma \protect\ref{LemmaP1} and (\protect\ref{eqP24})}
{Proof of Lemma 1 and (A.25)}}\label{app16}

The results in Lemma \ref{LemmaP1}, with the exception of (\ref
{eqP4}); and also result (\ref{eqP24}); will follow if we show that
for each $\ell\geq1$, some $C_1>0$ and all $C_2>0$,
%
%
\begin{eqnarray}
\label{eqP32}
P_k \Bigl\{\sup_{x\in\cI} \bigl|U_\ell(x)-
\ka_\ell f_X(x) \bigr|>n^{-C_1} \Bigr\} &=&O
\bigl(n^{-C_2} \bigr),
\\
\label{eqP33}
P_k \Bigl\{\max_{j=1,\ldots,n_k} \sup_{x\in\cI}
\bigl|U_{kj\ell}(x)-\ka_\ell f_X(x) \bigr|>n^{-C_1}
\Bigr\} &=&O \bigl(n^{-C_2} \bigr).
\end{eqnarray}
We shall derive (\ref{eqP33}); a proof of (\ref{eqP32}) is similar.

Markov's inequality can be used to prove that
%
%
\begin{equation}\label{eqP34}
\max_{j=1,\ldots,n_k} \sup_{x\in\cI} P_k
\bigl\{ \bigl|U_{kj\ell}(x)-\ka_\ell f_X(x)
\bigr|>n^{-C_1} \bigr\} =O \bigl(n^{-C_2} \bigr).
\end{equation}
It follows from (\ref{eq39})(c) that each $n_k$ is increasing no
faster than polynomially in $n$, and therefore, if we confine attention
to $x$ in a subset $\cI_n$, say, of $\cI$ that contains only $O(n^C)$
points for some $C>0$, we can place the maximum and supremum inside the
probability statement at (\ref{eqP34}), provided that $\cI$ is
replaced by $\cI_n$: for some $C_1>0$ and all $C_2>0$,
%
%
\begin{equation}\label{eqP35}
P_k \Bigl\{\max_{j=1,\ldots,n_k} \sup_{x\in\cI_n}
\bigl|U_{kj\ell}(x)-\ka_\ell f_X(x) \bigr|>n^{-C_1}
\Bigr\} =O \bigl(n^{-C_2} \bigr).
\end{equation}
The assumption, in (\ref{eq39})(b), that $K$ is compactly supported
and H\"older continuous, and the implication, in (\ref{eq37})(c),
that $f_X$ is also H\"older continuous, enable (\ref{eqP33}) to be
derived directly from (\ref{eqP35}) by taking $\cI_n$ to be a
sufficiently fine grid in $\cI$.

A proof of (\ref{eqP4}) in Lemma \ref{LemmaP1} is similar. To
illustrate the argument, we derive the following result part of (\ref
{eqP4}): for all $C_2,C_4>0$,
%
%
\begin{equation}\label{eqP36}
P_k \Bigl\{\sup_{x\in\cI} \bigl|\De(x)\bigr|>n^{C_4}
(mh)^{-1/2} \Bigr\} =O \bigl(n^{-C_2} \bigr).
\end{equation}
Using Markov's and Rosenthal's inequalities, we first obtain the result
when the supremum is outside the probability statement:
\[
\sup_{x\in\cI} P_k \bigl\{\bigl|\De(x)\bigr|>n^{C_4}
(mh)^{-1/2} \bigr\} =O \bigl(n^{-C_2} \bigr).
\]
Taking $\cI_n$ to contain only $O(n^C)$ points, for any fixed $C>0$,
we deduce that
\[
P_k \Bigl\{\sup_{x\in\cI_n} \bigl|\De(x)\bigr|>n^{C_4}
(mh)^{-1/2} \Bigr\} =O \bigl(n^{-C_2} \bigr),
\]
and taking $\cI_n$ to be a sufficiently fine grid in $\cI$ we obtain
(\ref{eqP36}).
\end{appendix}


\begin{supplement}
\stitle{Supplement to ``Unexpected properties of bandwidth
choice when smoothing discrete data for constructing a functional data
classifier''\\}
\slink[doi]{10.1214/13-AOS1158SUPP} 
\sdatatype{.pdf}
\sfilename{aos1158\_supp.pdf}
\sdescription{The supplementary file contains the proof of
Theorems \ref{Theorem2} and \ref{Theorem3}, as well as additional simulation results.}
\end{supplement}


\printaddresses


\begin{thebibliography}{30}

\bibitem[\protect\citeauthoryear{Araki et~al.}{2009}]{Araetal09}
\begin{barticle}[mr]
\bauthor{\bsnm{Araki},~\bfnm{Yuko}\binits{Y.}},
  \bauthor{\bsnm{Konishi},~\bfnm{Sadanori}\binits{S.}},
  \bauthor{\bsnm{Kawano},~\bfnm{Shuichi}\binits{S.}} \AND
  \bauthor{\bsnm{Matsui},~\bfnm{Hidetoshi}\binits{H.}}
(\byear{2009}).
\btitle{Functional logistic discrimination via regularized basis expansions}.
\bjournal{Comm. Statist. Theory Methods}
\bvolume{38}
\bpages{2944--2957}.
\bid{doi={10.1080/03610920902947246}, issn={0361-0926}, mr={2568196}}
\bptok{imsref}%
\end{barticle}
\endbibitem

\bibitem[\protect\citeauthoryear{Benhennia and Degras}{2011}]{BenDeg}
\begin{bmisc}[auto:STB|2013/09/19|12:14:10]
\bauthor{\bsnm{Benhennia},~\bfnm{K.}\binits{K.}} \AND
  \bauthor{\bsnm{Degras},~\bfnm{D.}\binits{D.}}
(\byear{2011}).
\bhowpublished{Local polynomial regression based on functional data.
  Unpublished manuscript. Available at \url{http://arxiv.org/pdf/1107.4058v1}}.
\bptok{imsref}%
\end{bmisc}
\endbibitem

\bibitem[\protect\citeauthoryear{Berlinet, Biau and Rouvi{\`e}re}{2008}]{BerBiaRou08}
\begin{barticle}[mr]
\bauthor{\bsnm{Berlinet},~\bfnm{Alain}\binits{A.}},
  \bauthor{\bsnm{Biau},~\bfnm{G{\'e}rard}\binits{G.}} \AND
  \bauthor{\bsnm{Rouvi{\`e}re},~\bfnm{Laurent}\binits{L.}}
(\byear{2008}).
\btitle{Functional supervised classification with wavelets}.
\bjournal{Ann. I.S.U.P.}
\bvolume{52}
\bpages{61--80}.
\bid{mr={2435041}}
\bptok{imsref}%
\end{barticle}
\endbibitem

\bibitem[\protect\citeauthoryear{Biau, Bunea and Wegkamp}{2005}]{BiaBunWeg05}
\begin{barticle}[mr]
\bauthor{\bsnm{Biau},~\bfnm{G{\'e}rard}\binits{G.}},
  \bauthor{\bsnm{Bunea},~\bfnm{Florentina}\binits{F.}} \AND
  \bauthor{\bsnm{Wegkamp},~\bfnm{Marten~H.}\binits{M.~H.}}
(\byear{2005}).
\btitle{Functional classification in {H}ilbert spaces}.
\bjournal{IEEE Trans. Inform. Theory}
\bvolume{51}
\bpages{2163--2172}.
\bid{doi={10.1109/TIT.2005.847705}, issn={0018-9448}, mr={2235289}}
\bptok{imsref}%
\end{barticle}
\endbibitem

\bibitem[\protect\citeauthoryear{Cardot, Degras and Josserand}{2013}]{CarDegJos}
\begin{barticle}[auto:STB|2013/09/19|12:14:10]
\bauthor{\bsnm{Cardot},~\bfnm{H.}\binits{H.}},
  \bauthor{\bsnm{Degras},~\bfnm{D.}\binits{D.}} \AND
  \bauthor{\bsnm{Josserand},~\bfnm{E.}\binits{E.}}
(\byear{2013}).
\btitle{Confidence bands for Horvitz--Thompson estimators using sampled
  noisy functional data}.
\bjournal{Bernoulli}
\bvolume{19}
\bpages{2067--2097}.
\bptok{imsref}%
\end{barticle}
\endbibitem

\bibitem[\protect\citeauthoryear{Cardot and Josserand}{2011}]{CarJos11}
\begin{barticle}[mr]
\bauthor{\bsnm{Cardot},~\bfnm{Herv{\'e}}\binits{H.}} \AND
  \bauthor{\bsnm{Josserand},~\bfnm{Etienne}\binits{E.}}
(\byear{2011}).
\btitle{Horvitz--{T}hompson estimators for functional data: Asymptotic
  confidence bands and optimal allocation for stratified sampling}.
\bjournal{Biometrika}
\bvolume{98}
\bpages{107--118}.
\bid{doi={10.1093/biomet/asq070}, issn={0006-3444}, mr={2804213}}
\bptok{imsref}%
\end{barticle}
\endbibitem

\bibitem[\protect\citeauthoryear{Carroll, Delaigle and Hall}{2013}]{CarDelHal}
\begin{bmisc}[auto:STB|2013/09/19|12:14:10]
\bauthor{\bsnm{Carroll},~\bfnm{R.~J.}\binits{R.~J.}},
  \bauthor{\bsnm{Delaigle},~\bfnm{A.}\binits{A.}} \AND
  \bauthor{\bsnm{Hall},~\bfnm{P.}\binits{P.}}
(\byear{2013}).
\bhowpublished{Supplement to ``Unexpected properties of bandwidth choice when
  smoothing discrete data for constructing a functional data classifier.''
  DOI:\doiurl{10.1214/13-AOS1158SUPP}}.
\bptok{imsref}%
\end{bmisc}
\endbibitem

\bibitem[\protect\citeauthoryear{Cuevas, Febrero and Fraiman}{2007}]{CueFebFra07}
\begin{barticle}[mr]
\bauthor{\bsnm{Cuevas},~\bfnm{Antonio}\binits{A.}},
  \bauthor{\bsnm{Febrero},~\bfnm{Manuel}\binits{M.}} \AND
  \bauthor{\bsnm{Fraiman},~\bfnm{Ricardo}\binits{R.}}
(\byear{2007}).
\btitle{Robust estimation and classification for functional data via
  projection-based depth notions}.
\bjournal{Comput. Statist.}
\bvolume{22}
\bpages{481--496}.
\bid{doi={10.1007/s00180-007-0053-0}, issn={0943-4062}, mr={2336349}}
\bptok{imsref}%
\end{barticle}
\endbibitem

\bibitem[\protect\citeauthoryear{Delaigle and Hall}{2012}]{DelHal12}
\begin{barticle}[mr]
\bauthor{\bsnm{Delaigle},~\bfnm{Aurore}\binits{A.}} \AND
  \bauthor{\bsnm{Hall},~\bfnm{Peter}\binits{P.}}
(\byear{2012}).
\btitle{Achieving near perfect classification for functional data}.
\bjournal{J. R. Stat. Soc. Ser. B Stat. Methodol.}
\bvolume{74}
\bpages{267--286}.
\bid{doi={10.1111/j.1467-9868.2011.01003.x}, issn={1369-7412}, mr={2899863}}
\bptok{imsref}%
\end{barticle}
\endbibitem

\bibitem[\protect\citeauthoryear{Delaigle, Hall and Bathia}{2012}]{DelHalBat12}
\begin{barticle}[mr]
\bauthor{\bsnm{Delaigle},~\bfnm{A.}\binits{A.}},
  \bauthor{\bsnm{Hall},~\bfnm{P.}\binits{P.}} \AND
  \bauthor{\bsnm{Bathia},~\bfnm{N.}\binits{N.}}
(\byear{2012}).
\btitle{Componentwise classification and clustering of functional data}.
\bjournal{Biometrika}
\bvolume{99}
\bpages{299--313}.
\bid{doi={10.1093/biomet/ass003}, issn={0006-3444}, mr={2931255}}
\bptok{imsref}%
\end{barticle}
\endbibitem

\bibitem[\protect\citeauthoryear{Epifanio}{2008}]{Epi08}
\begin{barticle}[mr]
\bauthor{\bsnm{Epifanio},~\bfnm{Irene}\binits{I.}}
(\byear{2008}).
\btitle{Shape descriptors for classification of functional data}.
\bjournal{Technometrics}
\bvolume{50}
\bpages{284--294}.
\bid{doi={10.1198/004017008000000154}, issn={0040-1706}, mr={2528652}}
\bptok{imsref}%
\end{barticle}
\endbibitem

\bibitem[\protect\citeauthoryear{Fan and Gijbels}{1996}]{FanGij96}
\begin{bbook}[mr]
\bauthor{\bsnm{Fan},~\bfnm{J.}\binits{J.}} \AND
  \bauthor{\bsnm{Gijbels},~\bfnm{I.}\binits{I.}}
(\byear{1996}).
\btitle{Local Polynomial Modelling and Its Applications}.
\bseries{Monographs on Statistics and Applied Probability}
\bvolume{66}.
\bpublisher{Chapman \& Hall}, \blocation{London}.
\bid{mr={1383587}}
\bptok{imsref}%
\end{bbook}
\endbibitem

\bibitem[\protect\citeauthoryear{Fromont and Tuleau}{2006}]{FroTul06}
\begin{bincollection}[mr]
\bauthor{\bsnm{Fromont},~\bfnm{Magalie}\binits{M.}} \AND
  \bauthor{\bsnm{Tuleau},~\bfnm{Christine}\binits{C.}}
(\byear{2006}).
\btitle{Functional classification with margin conditions}.
In \bbooktitle{Learning Theory---Proceedings of the 19th Annual Conference on Learning Theory,
Pittsburgh, 2006}
(\beditor{J. G. Carbonell} and \beditor{J. Siekmann}, eds.).
\bpublisher{Springer}, \blocation{New York}.
\bptok{imsref}%
\end{bincollection}
\endbibitem

\bibitem[\protect\citeauthoryear{Hall and Hosseini-Nasab}{2006}]{HalHos06}
\begin{barticle}[mr]
\bauthor{\bsnm{Hall},~\bfnm{Peter}\binits{P.}} \AND
  \bauthor{\bsnm{Hosseini-Nasab},~\bfnm{Mohammad}\binits{M.}}
(\byear{2006}).
\btitle{On properties of functional principal components analysis}.
\bjournal{J. R. Stat. Soc. Ser. B Stat. Methodol.}
\bvolume{68}
\bpages{109--126}.
\bid{doi={10.1111/j.1467-9868.2005.00535.x}, issn={1369-7412}, mr={2212577}}
\bptok{imsref}%
\end{barticle}
\endbibitem

\bibitem[\protect\citeauthoryear{Hall and Hosseini-Nasab}{2009}]{HalHos09}
\begin{barticle}[mr]
\bauthor{\bsnm{Hall},~\bfnm{Peter}\binits{P.}} \AND
  \bauthor{\bsnm{Hosseini-Nasab},~\bfnm{Mohammad}\binits{M.}}
(\byear{2009}).
\btitle{Theory for high-order bounds in functional principal components
  analysis}.
\bjournal{Math. Proc. Cambridge Philos. Soc.}
\bvolume{146}
\bpages{225--256}.
\bid{doi={10.1017/S0305004108001850}, issn={0305-0041}, mr={2461880}}
\bptok{imsref}%
\end{barticle}
\endbibitem

\bibitem[\protect\citeauthoryear{Hall and Kang}{2005}]{HalKan05}
\begin{barticle}[mr]
\bauthor{\bsnm{Hall},~\bfnm{Peter}\binits{P.}} \AND
  \bauthor{\bsnm{Kang},~\bfnm{Kee-Hoon}\binits{K.-H.}}
(\byear{2005}).
\btitle{Bandwidth choice for nonparametric classification}.
\bjournal{Ann. Statist.}
\bvolume{33}
\bpages{284--306}.
\bid{doi={10.1214/009053604000000959}, issn={0090-5364}, mr={2157804}}
\bptok{imsref}%
\end{barticle}
\endbibitem

\bibitem[\protect\citeauthoryear{Hall and Van~Keilegom}{2007}]{HalVan07}
\begin{barticle}[mr]
\bauthor{\bsnm{Hall},~\bfnm{Peter}\binits{P.}} \AND
  \bauthor{\bsnm{Van~Keilegom},~\bfnm{Ingrid}\binits{I.}}
(\byear{2007}).
\btitle{Two-sample tests in functional data analysis starting from discrete
  data}.
\bjournal{Statist. Sinica}
\bvolume{17}
\bpages{1511--1531}.
\bid{issn={1017-0405}, mr={2413533}}
\bptok{imsref}%
\end{barticle}
\endbibitem

\bibitem[\protect\citeauthoryear{Leng and M{\"{u}}ller}{2006}]{LenMul06}
\begin{barticle}[pbm]
\bauthor{\bsnm{Leng},~\bfnm{Xiaoyan}\binits{X.}} \AND
  \bauthor{\bsnm{M{\"{u}}ller},~\bfnm{Hans-Georg}\binits{H.-G.}}
(\byear{2006}).
\btitle{Classification using functional data analysis for temporal gene
  expression data}.
\bjournal{Bioinformatics}
\bvolume{22}
\bpages{68--76}.
\bid{doi={10.1093/bioinformatics/bti742}, issn={1367-4803}, pii={bti742},
  pmid={16257986}}
\bptok{imsref}%
\end{barticle}
\endbibitem

\bibitem[\protect\citeauthoryear{Li and Hsing}{2010a}]{LiHsi10N1}
\begin{barticle}[mr]
\bauthor{\bsnm{Li},~\bfnm{Yehua}\binits{Y.}} \AND
  \bauthor{\bsnm{Hsing},~\bfnm{Tailen}\binits{T.}}
(\byear{2010}a).
\btitle{Deciding the dimension of effective dimension reduction space for
  functional and high-dimensional data}.
\bjournal{Ann. Statist.}
\bvolume{38}
\bpages{3028--3062}.
\bid{doi={10.1214/10-AOS816}, issn={0090-5364}, mr={2722463}}
\bptok{imsref}%
\end{barticle}
\endbibitem

\bibitem[\protect\citeauthoryear{Li and Hsing}{2010b}]{LiHsi10N2}
\begin{barticle}[mr]
\bauthor{\bsnm{Li},~\bfnm{Yehua}\binits{Y.}} \AND
  \bauthor{\bsnm{Hsing},~\bfnm{Tailen}\binits{T.}}
(\byear{2010}b).
\btitle{Uniform convergence rates for nonparametric regression and principal
  component analysis in functional/longitudinal data}.
\bjournal{Ann. Statist.}
\bvolume{38}
\bpages{3321--3351}.
\bid{doi={10.1214/10-AOS813}, issn={0090-5364}, mr={2766854}}
\bptok{imsref}%
\end{barticle}
\endbibitem

\bibitem[\protect\citeauthoryear{L{\'o}pez-Pintado and Romo}{2006}]{LopRom06}
\begin{bincollection}[mr]
\bauthor{\bsnm{L{\'o}pez-Pintado},~\bfnm{Sara}\binits{S.}} \AND
  \bauthor{\bsnm{Romo},~\bfnm{Juan}\binits{J.}}
(\byear{2006}).
\btitle{Depth-based classification for functional data}.
In \bbooktitle{Data Depth: Robust Multivariate Analysis, Computational Geometry
  and Applications}.
\bseries{DIMACS Series in Discrete Mathematics and Theoretical Computer Science}
\bvolume{72}
\bpages{103--119}.
\bpublisher{Amer. Math. Soc.}, \blocation{Providence, RI}.
\bid{mr={2343116}}
\bptok{imsref}%
\end{bincollection}
\endbibitem

\bibitem[\protect\citeauthoryear{Manning, Raghavan and Sch{\"u}tze}{2008}]{ManRagSch08}
\begin{bbook}[auto:STB|2013/09/19|12:14:10]
\bauthor{\bsnm{Manning},~\bfnm{C.~D.}\binits{C.~D.}},
  \bauthor{\bsnm{Raghavan},~\bfnm{P.}\binits{P.}} \AND
  \bauthor{\bsnm{Sch{\"u}tze},~\bfnm{H.}\binits{H.}}
(\byear{2008}).
\btitle{Introduction to Information Retrival}.
\bpublisher{Cambridge Univ. Press}, \blocation{Cambridge}.
\bptok{imsref}%
\end{bbook}
\endbibitem

\bibitem[\protect\citeauthoryear{Panaretos, Kraus and Maddocks}{2010}]{PanKraMad10}
\begin{barticle}[mr]
\bauthor{\bsnm{Panaretos},~\bfnm{Victor~M.}\binits{V.~M.}},
  \bauthor{\bsnm{Kraus},~\bfnm{David}\binits{D.}} \AND
  \bauthor{\bsnm{Maddocks},~\bfnm{John~H.}\binits{J.~H.}}
(\byear{2010}).
\btitle{Second-order comparison of {G}aussian random functions and the geometry
  of {DNA} minicircles}.
\bjournal{J. Amer. Statist. Assoc.}
\bvolume{105}
\bpages{670--682}.
\bid{doi={10.1198/jasa.2010.tm09239}, issn={0162-1459}, mr={2724851}}
\bptok{imsref}%
\end{barticle}
\endbibitem

\bibitem[\protect\citeauthoryear{Ramsay and Silverman}{2005}]{RamSil05}
\begin{bbook}[mr]
\bauthor{\bsnm{Ramsay},~\bfnm{J.~O.}\binits{J.~O.}} \AND
  \bauthor{\bsnm{Silverman},~\bfnm{B.~W.}\binits{B.~W.}}
(\byear{2005}).
\btitle{Functional Data Analysis},
\bedition{2nd} ed.
\bpublisher{Springer}, \blocation{New York}.
\bid{mr={2168993}}
\bptok{imsref}%
\end{bbook}
\endbibitem

\bibitem[\protect\citeauthoryear{Rossi and Villa}{2006}]{RosVil06}
\begin{barticle}[auto:STB|2013/09/19|12:14:10]
\bauthor{\bsnm{Rossi},~\bfnm{F.}\binits{F.}} \AND
  \bauthor{\bsnm{Villa},~\bfnm{N.}\binits{N.}}
(\byear{2006}).
\btitle{Support vector machine for functional data classification}.
\bjournal{Neurocomputing}
\bvolume{69}
\bpages{730--742}.
\bptok{imsref}%
\end{barticle}
\endbibitem

\bibitem[\protect\citeauthoryear{Ruppert, Sheather and Wand}{1995}]{RupSheWan95}
\begin{barticle}[mr]
\bauthor{\bsnm{Ruppert},~\bfnm{D.}\binits{D.}},
  \bauthor{\bsnm{Sheather},~\bfnm{S.~J.}\binits{S.~J.}} \AND
  \bauthor{\bsnm{Wand},~\bfnm{M.~P.}\binits{M.~P.}}
(\byear{1995}).
\btitle{An effective bandwidth selector for local least squares regression}.
\bjournal{J. Amer. Statist. Assoc.}
\bvolume{90}
\bpages{1257--1270}.
\bid{issn={0162-1459}, mr={1379468}}
\bptok{imsref}%
\end{barticle}
\endbibitem

\bibitem[\protect\citeauthoryear{Vilar and P{\'e}rtega}{2004}]{VilPer04}
\begin{barticle}[mr]
\bauthor{\bsnm{Vilar},~\bfnm{Jos{\'e}~A.}\binits{J.~A.}} \AND
  \bauthor{\bsnm{P{\'e}rtega},~\bfnm{Sonia}\binits{S.}}
(\byear{2004}).
\btitle{Discriminant and cluster analysis for {G}aussian stationary processes:
  Local linear fitting approach}.
\bjournal{J. Nonparametr. Stat.}
\bvolume{16}
\bpages{443--462}.
\bid{doi={10.1080/10485250410001656453}, issn={1048-5252}, mr={2073035}}
\bptok{imsref}%
\end{barticle}
\endbibitem

\bibitem[\protect\citeauthoryear{Wang, Ray and Mallick}{2007}]{WanRayMal07}
\begin{barticle}[mr]
\bauthor{\bsnm{Wang},~\bfnm{Xiaohui}\binits{X.}},
  \bauthor{\bsnm{Ray},~\bfnm{Shubhankar}\binits{S.}} \AND
  \bauthor{\bsnm{Mallick},~\bfnm{Bani~K.}\binits{B.~K.}}
(\byear{2007}).
\btitle{Bayesian curve classification using wavelets}.
\bjournal{J.~Amer. Statist. Assoc.}
\bvolume{102}
\bpages{962--973}.
\bid{doi={10.1198/016214507000000455}, issn={0162-1459}, mr={2354408}}
\bptok{imsref}%
\end{barticle}
\endbibitem

\bibitem[\protect\citeauthoryear{Wu and M{\"u}ller}{2011}]{WuMul11}
\begin{barticle}[mr]
\bauthor{\bsnm{Wu},~\bfnm{Shuang}\binits{S.}} \AND
  \bauthor{\bsnm{M{\"u}ller},~\bfnm{Hans-Georg}\binits{H.-G.}}
(\byear{2011}).
\btitle{Response-adaptive regression for longitudinal data}.
\bjournal{Biometrics}
\bvolume{67}
\bpages{852--860}.
\bid{doi={10.1111/j.1541-0420.2010.01518.x}, issn={0006-341X}, mr={2829259}}
\bptok{imsref}%
\end{barticle}
\endbibitem

\end{thebibliography}
\end{document}